\input ./style/arxiv-general.cfg
\documentclass[MSNbibl,number,citesort,seceqn,dvips]{arxbj}
\makeatletter
   \@ifpackageloaded{graphicx}{}{\usepackage{graphicx}}
\makeatother
\usepackage{mathbh,dcolumn}

\volume{23}
\issue{1}
\pubyear{2017}
\firstpage{191}
\lastpage{218}
\doi{10.3150/15-BEJ741}
\docsubty{FLA}

\makeatletter
\newcolumntype{d}[1]{D{.}{.}{#1}}

\newcommand{\rrvert}{\vert}

\newcommand{\llvert}{\vert}
\def\overset{\stackrel}
\def\mathbbm{\mathbh}
\newcommand{\eqref}[1]{(\ref{#1})}

\def\bigintsss{\int}
\newtheorem{thmm}{Theorem}[section]
\newremark{rem}{Remark}[section]
\newremark{remarks}{Remarks}

\newproclaim{definition}{Definition}
\newremark{example}{Example}

\newproclaim{assumption}{Assumption}

\newcommand{\Prob}{\mathbb{P} }

\newtheorem{Prop}{Proposition}[section]
\newtheorem{Lemma}{Lemma}[section]

\newtheorem{Corollary}{Corollary}[section]

\makeatother

\begin{document}
\begin{frontmatter}

\title{Bounds for the normal approximation of the maximum likelihood estimator}
\runtitle{Distance to normal for MLE}

\begin{aug}
\author[A]{\inits{A.}\fnms{Andreas} \snm{Anastasiou}\corref{}\thanksref
{A,e1}\ead
[label=e1,mark]{anastasi@stats.ox.ac.uk}}
\and
\author[A]{\inits{G.}\fnms{Gesine} \snm{Reinert}\thanksref
{A,e2}\ead
[label=e2,mark]{reinert@stats.ox.ac.uk}}
\address[A]{Department of Statistics, University of Oxford, 1 South
Parks Road,
Oxford,
OX1 3TG,
UK.\\
\printead{e1,e2}}
\end{aug}

%
\received{\smonth{12} \syear{2014}}
%
\revised{\smonth{6} \syear{2015}}


\begin{abstract}
While the asymptotic normality of the maximum likelihood estimator
under regularity conditions is long established, this paper derives
explicit bounds for the bounded Wasserstein distance between the
distribution of the maximum likelihood estimator (MLE) and the normal
distribution. For this task, we employ Stein's method. We focus on
independent and identically distributed random variables, covering both
discrete and continuous distributions as well as exponential and
non-exponential families. In particular, a closed form expression of
the MLE is not required. We also use a perturbation method to treat
cases where the MLE has positive probability of being on the boundary
of the parameter space.
\end{abstract}

%
\begin{keyword}
\kwd{maximum likelihood estimator}
\kwd{normal approximation}
\kwd{Stein's method}
\end{keyword}
\end{frontmatter}

\section{Introduction}
\label{sec:intro} This paper assesses the bounded Wasserstein distance
between the distribution of the maximum likelihood estimator (MLE) and
the normal distribution. We concentrate on independent and identically
distributed (i.i.d.) random variables, with the case that the random
variables follow an exponential family distribution as an example. We
also explain how a perturbation of both the parameter and the data can
be useful in specific situations. The treatment includes situations
where the MLE has positive probability to be on the boundary of the
parameter space. The paper also covers cases where there is not an
analytic form for the MLE.

Here is the notation which is used throughout the paper. First of all,
$\theta$ denotes a scalar unknown parameter found in a parametric
statistical model. Let $\theta_0$ be the true (still unknown) value of
the parameter $\theta$ and let $\Theta\subset\mathbb{R}$ denote the
parameter space, while $\mathbf{X} = (X_1, X_2, \ldots, X_n)$ is
the random sample of $n$ i.i.d. random variables with joint density
function $f(\mathbf{x}|\theta)$. For $X_i = x_i$ being some
observed values, the likelihood function is $L(\theta; \mathbf{x})
= f(\mathbf{x}|\theta)$. Its natural logarithm, called the
log-likelihood function is denoted by $l(\theta;\mathbf{x})$.
Having a fixed set of data and a defined statistical model, a maximum
likelihood estimate is a value of the parameter which maximises the
likelihood function. Derivatives of the log-likelihood function, with
respect to $\theta$, are denoted by $l'(\theta;\mathbf
{x}),l''(\theta;\mathbf{x}),\ldots, l^{(j)}(\theta;\mathbf
{x})$, for $j$ any integer greater than 2. For many models, the MLE
exists and it is also unique, in which case it is denoted by $\hat
{\theta}_n(\mathbf{X})$; this is known as the ``regular'' case.
However, uniqueness or even existence of the MLE is not always secured.
Unless otherwise specified, we make the following assumptions:
\begin{enumerate}[(iii)]
\item[(i)] The log-likelihood function $l(\theta;\mathbf{x})$ is a
twice continuously differentiable function with respect to $\theta$ and
the parameter varies in an open interval $(a,b)$, where $a, b
\in\mathbb{R}\cup \{-\infty,\infty \}$ and $a < b$.
\item[(ii)] $\lim_{\theta\to a, b} l(\theta;\mathbf{x}) =
-\infty$,
\item[(iii)] $l''(\theta;\mathbf{x}) < 0$ at every point $\theta
\in(a,b)$ for which $l'(\theta;\mathbf{x}) = 0$.
\end{enumerate}

Under the assumptions (i)--(iii) above, the MLE exists and it is
unique (Makelainen
\textit{et~al.}~\cite{Makelainen}).
Following now Casella and Berger \cite{Casella}, unless otherwise
stated we also make the
following assumptions:
\begin{enumerate}[(R1)]
\item[(R1)] the parameter is identifiable, which means that if $\theta
\neq\theta'$, then $\exists x: f(x|\theta) \not\equiv f(x|\theta')$;
\item[(R2)]
the density $f(x|\theta)$ is three times differentiable with respect to
$\theta$, the third derivative is continuous in $\theta$ and
$\bigintsss f(x|\theta) \,\mathrm{d}x$ can be differentiated three times
under the integral sign;
\item[(R3)] for any $\theta_0 \in\Theta$ and for $\mathbb{X}$ denoting
the support of $f(x|\theta)$, there exists a positive number
$\varepsilon$
and a function $M(x)$ (both of which may depend on $\theta_0$) such that
\[
\biggl\llvert \frac{\mathrm{d}^3}{\mathrm{d}\theta^3}\log
f(x|\theta )\biggr\rrvert \leq M(x)\qquad
\forall x \in\mathbb{X}, \theta_0 - \varepsilon< \theta<
\theta_0 + \varepsilon,
\]
with ${\mathrm{E}}_{\theta_0}[M(X)] < \infty$;
\item[(R4)] $i(\theta_0) \neq0$, where $i(\theta)$ is the expected
Fisher Information for one random variable.
\end{enumerate}
The requirement (R2) that $\bigintsss f(x|\theta) \,\mathrm{d}x$ can be
differentiated three times under the integral sign is usually
substituted in the literature by the assumption that integration of
$f(x|\theta)$ over $x$ and differentiation with respect to $\theta$ are
three times interchangeable, so that $\bigintsss_{\mathbb{R}}^{}\frac
{\mathrm{d}^j}{\mathrm{d}\theta^j} f(x|\theta) \,\mathrm{d}x = \frac
{\mathrm{d}^j}{\mathrm{d}\theta^j}\bigintsss_{\mathbb{R}}^{}
f(x|\theta
) \,\mathrm{d}x = 0, j \in \{1, 2, 3 \}$. This
condition ensures that if the expressions exist, then ${\mathrm
{E}}_{\theta
}[l'(\theta;\mathbf{X})] = 0$ and $\operatorname{Var}_{\theta
}[l'(\theta
;\mathbf{X})] = n i(\theta)$. In addition, it is obvious from (R3)
that $ \{\theta: |\theta- \theta_0| < \varepsilon \}
\subset\Theta$ is required. The motivation of the work presented in
this paper are the results given in Theorem~\ref{MLEt}. The efficiency
and asymptotic normality of the MLE have first been discussed in Fisher
\cite
{Fisher}. Here we present the i.i.d. case; see Hoadley \cite{Hoadley}
for the
case of independent but not identically distributed random variables.

%
\begin{thmm} [(Casella and Berger \cite{Casella}, page~472)]
\label{MLEt}
Let $X_1, X_2, \ldots, X_n$ be i.i.d. random variables with probability
density (or mass) function $f(x_i|\theta)$, where $\theta$ is the
scalar parameter. Assume that the MLE exists and it is unique and
\textup{(R1)--(R4)} are satisfied. Then for $Z \sim{\mathrm{N}}(0,1)$,
%
\begin{equation}
\label{scoredistribution} \mathrm{(a)}\quad \frac{1}{\sqrt{n}}l'(\theta_0;
\mathbf{X}) \mathop{\longrightarrow}^{\mathrm{d}}_{n \to\infty} \sqrt{i(
\theta_0)}Z,\qquad \mathrm{(b)}\quad \sqrt{n i(\theta _0)}\bigl(\hat{
\theta}_n(\mathbf{X}) - \theta_0\bigr)
\mathop{\longrightarrow}^{\mathrm{d}}_{n \to \infty} Z.
\end{equation}
\end{thmm}
Theorem~\ref{MLEt} gives only a qualitative result as $n \rightarrow
\infty$, but in approximations the sample size, $n$, is always finite
and it is not clear when $n$ is ``large enough'' for the limiting
behaviour to be a good approximation to the finite-$n$ behaviour. The
rate of convergence may also depend on the true parameter $\theta_0$.
Hence, it is of interest to obtain explicit bounds for a distributional
distance related to (a) and (b) in \eqref{scoredistribution}. These
bounds are given in Proposition~\ref{propscore} and Theorem~\ref
{Theoremnoncan}, respectively. The tools we use are mainly Taylor
expansions, conditional expectations, a perturbation method and a
result from Stein's method as given in Lemma~\ref{Gesinetheorem}.
Bounds are also derived in Geyer \cite{Geyer}, using the framework of locally
asymptotically mixed normal (LAMN) models, but these bounds are of
asymptotic nature.

As distance, we mainly use the bounded Wasserstein distance. If $F, G$
are two random variables with values in $\mathbb{R}$ and $H$ is a class
of separating functions, then a Zolotarev-type distance between the
laws of $F$ and $G$, induced by $H$, is given by the quantity
%
\begin{equation}
\label{nourdingeneral} d_{H}(F,G) = \sup \bigl\{\bigl|{\mathrm{E}}\bigl[h(F)\bigr]
- {\mathrm{E}}\bigl[h(G)\bigr]\bigr|: h \in H \bigr\}.
\end{equation}
From now on, $\|\cdot\|$ denotes the supremum norm ($\|\cdot\|_{\infty}$) and
%
\begin{equation}
\label{classfunctions} H = \bigl\{ h:\mathbb{R}\rightarrow\mathbb{R}:\|h\|_{{\mathrm{ Lip}}} +
\|h\| \leq1 \bigr\},
\end{equation}
where
\[
\|h\|_{{\mathrm{ Lip}}} = \mathop{\mathop{\sup}_{x,y \in
\mathbb{R}}}_{x
\neq y}\frac{|h(x) - h(y)|}{|x-y|}.
\]
Using Rademacher's
theorem, since $\|h\|_{{\mathrm{ Lip}}} \leq1$, then $h$ is differentiable
almost everywhere, with $h'$ denoting its derivative.

Using this class of test functions, \eqref{nourdingeneral} gives the
bounded Wasserstein (or Fortet--Mourier) distance between two random
variables $F$ and $G$, denoted from now on by
%
\begin{equation}
\label{general_bounded_Wasserstein} d_{bW}(F,G) = \sup \bigl\{\bigl|{\mathrm{E}}\bigl[h(F)\bigr]
- {\mathrm{E}}\bigl[h(G)\bigr]\bigr|: h \in H \bigr\},
\end{equation}
with $H$ as in \eqref{classfunctions}; see, for example, Nourdin and
Peccati \cite
{Nourdin}. Rachev \cite{Rachev} also gives a connection to the
Kantorovich--Rubinstein problem. To obtain such bounds, we use the
following lemma from Reinert \cite{Gesinepaper} which is based on Stein's
method (Stein \cite{Stein1972}).

\begin{Lemma}
\label{Gesinetheorem}
Let $Y_1, Y_2, \ldots, Y_n$ be independent random variables with
${\mathrm{E}}(Y_i) = 0$, $\operatorname{Var}(Y_i) = \sigma^2 > 0$
and ${\mathrm{E}}\llvert Y_i\rrvert ^3 < \infty$. Let $W=\frac
{1}{\sqrt{n}}\sum_{i=1}^{n}
Y_i$ and $K \sim{\mathrm{N}}(0,\sigma^2)$. Then for any function $h
\in H$,
with $H$ given in \eqref{classfunctions}
%
\begin{equation}
\label{equationLemma} d_{bW}(W,K) \leq\frac{1}{\sqrt{n}} \biggl(2 +
\frac{1}{\sigma
^3} \bigl[{\mathrm{E}}\llvert Y_1\rrvert ^3
\bigr] \biggr).
\end{equation}
\end{Lemma}
Using $Y_i = l'(\theta_0;X_i)$, we see that \eqref
{equationLemma} is closely related to (a) in \eqref
{scoredistribution}. For a bound of (b), we employ Taylor
expansion.

The paper is organised as follows. Section~\ref{sec:generalc} gives an
upper bound on the distributional distance between the distribution of
the MLE and the normal distribution in the case of i.i.d. random
variables. In Section~\ref{sec:one_parameter}, the results are applied
to the class of one-parameter exponential family distributions. In
Section~\ref{sec:perturbationc}, we use a perturbation to treat the
special case of having a random vector from a distribution where the
parameter space is not an open interval and there is positive
probability of the MLE to lie on the boundary of the parameter space.
An example is the Poisson distribution with mean $\theta\in[0,\infty
)$; the MLE could take on the value zero with positive probability, but
the log-likelihood function is not differentiable at zero. In
Section~\ref{sec:Implicit}, we obtain an upper bound on the Mean
Squared Error
of the MLE. We use this bound in order to get an upper bound on the
distributional distance to the normal distribution, even when no
analytic expression of the MLE is available. We assess the quality of
our results through a simulation-based study related to the Beta
distribution. The R-code for the simulations and the simulation output
are available at the Oxford University Research Archive (ORA). The DOI
is: \doiurl{10.5287/bodleian:s4655h876}.

\section{Bounds on the distance to normal for the MLE}
\label{sec:generalc}
In this section, we briefly relate the Kolmogorov and the bounded
Wasserstein distance and we give upper bounds on the distributional
distance between the distribution of the MLE and the normal
distribution in terms of the bounded Wasserstein distance.
\subsection{The bounded Wasserstein and the Kolmogorov distance}
For $Z \sim{\mathrm{N}}(0,1)$, the aim is to bound
%
\begin{equation}
\label{interest} d_{bW} \bigl(\sqrt{n i(\theta_0)}\bigl(
\hat{\theta}_n(\mathbf{X}) - \theta_0\bigr),Z
\bigr),
\end{equation}
with $d_{bW}(\cdot, \cdot)$ as defined in \eqref
{general_bounded_Wasserstein}.
Using $H =  \{{\mathbbm{1}_{[\cdot\leq x]}}, x \in\mathbb
{R} \}$ as the class of functions in \eqref{nourdingeneral},
yields the Kolmogorov distance,
\[
d_K \bigl(\sqrt{n i(\theta_0)}\bigl(\hat{
\theta}_n(\mathbf{X}) - \theta_0\bigr),Z \bigr).
\]
The next proposition links these two distances.

%
\begin{Prop}
\label{Proposition_Kolmogorov}
If $G$ is any real-valued random variable and $Z \sim \mathrm{N}(0,1)$, then
\[
d_{K}(G,Z) \leq2\sqrt{d_{bW}(G,Z)}.
\]
\end{Prop}
\begin{pf}
The proof of this proposition follows the proof of Theorem~3.3
of Chen
\textit{et~al.} \cite{Chen_book}, page~48. Let $z \in\mathbb{R}$ and for
$\alpha= \sqrt{d_{bW}(G,Z)}(2\pi)^{{1}/{4}}$, $z \in\mathbb{R}$,
let
\[
h_{\alpha}(w) = %
\cases{ 1, & \quad $\mbox{if $w \leq z$}$,
\vspace*{2pt}
\cr
1 + \displaystyle\frac{z-w}{\alpha}, &\quad  $\mbox{if $z<w\leq z+\alpha$}$,\vspace
*{2pt}
\cr
0, &\quad  $\mbox{if $w > z+ \alpha$}$ } %
\]
so that $h_{\alpha}$ is bounded Lipschitz with $\|h_{\alpha}\| \leq1$
and $\|h'_{\alpha}\| \leq\frac{1}{\alpha}$. By the triangle inequality,
\begin{eqnarray*}
\Prob (G \leq z ) - \Prob (Z\leq z ) &\leq& {\mathrm{E}}
\bigl[h_{\alpha}(G) \bigr] - {\mathrm{E}} \bigl[h_{\alpha}(Z) \bigr] +
{\mathrm{E}} \bigl[h_{\alpha}(Z) \bigr] - \Prob (Z \leq z )
\\
&\leq&\frac{d_{bW}(G,Z)}{\alpha} + \Prob (z \leq Z \leq z+\alpha )
\\
& \leq&\frac{d_{bW}(G,Z)}{\alpha} + \frac{\alpha}{\sqrt
{2\pi
}} \leq2\sqrt{d_{bW}(G,Z)}.
\end{eqnarray*}
Similarly $\Prob (G \leq z ) - \Prob (Z\leq z )
\geq
-2\sqrt{d_{bW}(G,Z)}$, which completes the proof.
\end{pf}

The Kolmogorov distance relates directly to exact conservative
confidence intervals. Our results on the bounded Wasserstein distance
and Proposition~\ref{Proposition_Kolmogorov} give that
\[
d_K \bigl(\sqrt{n i(\theta_0)}\bigl(\hat{
\theta}_n(\mathbf{X}) - \theta_0\bigr),Z \bigr)
\leq2\sqrt{B_{bW}} =: B_K,
\]
where $B_{bW}$ denotes the bound for the bounded Wasserstein distance
from Proposition~\ref{Proposition_Kolmogorov}. Therefore, for $y \in
\mathbb{R}$:
%
\begin{eqnarray}
\label{conf.int.1}
&&\bigl\llvert \Prob \bigl(\sqrt{n i(
\theta_0)} \bigl(\hat {\theta }_n(\mathbf{X}) -
\theta_0 \bigr)\leq y \bigr) - \Prob (Z\leq y )\bigr\rrvert \leq
B_K
\nonumber
\\[-8pt]
\\[-8pt]
\nonumber
&&\quad\Leftrightarrow\quad -B_K \leq\Prob \bigl(\sqrt{n i(
\theta_0)} \bigl(\hat {\theta}_n(\mathbf{X}) -
\theta_0 \bigr)\leq y \bigr) - \Prob (Z\leq y ) \leq
B_K.
\end{eqnarray}
For $\Phi^{-1}(\cdot)$ the quantile function for the standard normal
distribution, applying \eqref{conf.int.1} to $y = \Phi^{-1}
(\frac
{\alpha}{2} - B_K )$ and to $y = \Phi^{-1} (1 - \frac
{\alpha
}{2} + B_K )$ yields
\[
\Prob \biggl(\Phi^{-1} \biggl(\frac{\alpha}{2} -
B_K \biggr) \leq \sqrt{n i(\theta_0)} \bigl(\hat{
\theta}_n(\mathbf{X}) - \theta _0 \bigr)\leq
\Phi^{-1} \biggl(1-\frac{\alpha}{2}+B_K \biggr) \biggr)
\geq1 - \alpha.
\]
Hence, if the expected Fisher Information number for one random
variable, $i(\theta_0)$, is known, then
\[
\biggl(\hat{\theta}_n(\mathbf{X}) -
\frac{\Phi
^{-1}
(1-{\alpha}/{2}+B_K )}{\sqrt{n i(\theta_0)}}, \hat {\theta }_n(\mathbf{X}) -
\frac{\Phi^{-1} ({\alpha}/{2}-B_K
)}{\sqrt{n i(\theta_0)}} \biggr)
\]
is a conservative $100(1-\alpha)\%$ confidence interval for $\theta_0$.
\subsection{Bounds in terms of the bounded Wasserstein distance}
The bounded Wasserstein distance links in well with Stein's method
because the Lipschitz test functions are differentiable
almost everywhere. From now on, $\frac{\mathrm{d}}{\mathrm
{d}\theta}\operatorname{ log}f(X_1|\theta_0) := \frac{\mathrm
{d}}{\mathrm{d}\theta
}\operatorname{ log}f(X_1|\theta) |_{{\theta= \theta_0}}$. The next
two results provide a bound for (a) and (b) in \eqref
{scoredistribution}, respectively.

%
\begin{Prop}
\label{propscore}
Suppose $X_1, X_2, \ldots, X_n$ are i.i.d. random variables with
density or frequency function $f(x_i|\theta)$. Assume that \textup{(R1)--(R4)}
are satisfied, $Z \sim \mathrm{N}(0,1)$ and ${\mathrm{E}}\llvert \frac{\mathrm
{d}}{\mathrm{d}\theta}\operatorname{ log}f(X_1|\theta_0)\rrvert ^3$
exists. Then
for $h:\mathbb{R}\rightarrow\mathbb{R}$, such that $h$ is absolutely
continuous and bounded
%
\begin{eqnarray}
\label{equationpropA11} \biggl\llvert {\mathrm{E}} \biggl[h \biggl(\frac{l'(\theta_0;\mathbf
{X})}{\sqrt{n
i(\theta_0)}}
\biggr) \biggr] - {\mathrm{E}}\bigl[h(Z)\bigr]\biggr\rrvert \leq\frac
{\|h'\|
}{\sqrt{n}}
\biggl(2 + \frac{1}{[i(\theta_0)]^{{3}/{2}}} \biggl[{\mathrm{E}}\biggl\llvert \frac{\mathrm{d}}{\mathrm{d}\theta
}
\operatorname{ log}f(X_1|\theta _0)\biggr\rrvert
^3 \biggr] \biggr).
\end{eqnarray}
In particular,
%
\begin{eqnarray}
\label{equationpropA1} d_{bW} \biggl(\frac{l'(\theta_0;\mathbf{X})}{\sqrt{n i(\theta_0)}}, Z \biggr) \leq
\frac{1}{\sqrt{n}} \biggl(2 + \frac{1}{[i(\theta
_0)]^{{3}/{2}}} \biggl[{\mathrm{E}}\biggl\llvert
\frac{\mathrm{d}}{\mathrm
{d}\theta}\operatorname{ log}f(X_1|\theta_0)
\biggr\rrvert ^3 \biggr] \biggr).
\end{eqnarray}
\end{Prop}
\begin{pf}
Let
\[
Y_i = Y_i(X_i;\theta_0) = \biggl(\frac{\mathrm{d}}{\mathrm
{d}\theta
}\operatorname{ log}f(X_i|\theta_0)\biggr)\Big/{\sqrt{i(\theta_0)}},\qquad
i=1,2,\ldots,n,
\]
which
are i.i.d. random variables as $X_1, X_2, \ldots, X_n$ are i.i.d. The
regularity conditions (R1)--(R4)\vspace*{1pt} ensure that ${\mathrm{E}}_{\theta
_0}[Y_i] =
0$ and $\operatorname{Var}_{\theta_0}[Y_i] = 1$. Then letting\vspace*{1pt}
$W=W(\mathbf
{X};\theta_0) = \frac{1}{\sqrt{n}}\sum_{i=1}^{n}Y_i = \frac
{l'(\theta
_0;\mathbf{X})}{\sqrt{n i(\theta_0)}}$, gives that ${\mathrm
{E}}_{\theta
_0}[W] = 0$ and $\operatorname{Var}_{\theta_0}[W] = 1$. Applying
Lemma~\ref
{Gesinetheorem} to $K = Z \sim \mathrm{N}(0,1)$ yields the result.
\end{pf}
%

%
\begin{thmm}
\label{Theoremnoncan}
Let $X_1, X_2, \ldots, X_n$ be i.i.d. random variables with density or
frequency function $f(x_i|\theta)$ such that the regularity conditions
\textup{(R1)--(R4)} are satisfied and that the MLE, $\hat{\theta}_n(\mathbf
{X})$, exists and it is unique. Assume that ${\mathrm{E}}\llvert \frac
{\mathrm
{d}}{\mathrm{d}\theta}\operatorname{ log}f(X_1|\theta_0)\rrvert ^3
< \infty$ and
that ${\mathrm{E}} (\hat{\theta}_n(\mathbf{X}) - \theta
_0 )^4 <
\infty$. Let $0 < \varepsilon= \varepsilon(\theta_0)$ be such that
$(\theta
_0 - \varepsilon, \theta_0 + \varepsilon) \subset\Theta$ as in
\textup{(R3)} and let
$Z \sim{\mathrm{N}}(0,1)$. Then
%
\begin{eqnarray}
\label{boundTHEOREM}
&& d_{bW} \bigl(\sqrt{n i(\theta_0)}
\bigl(\hat{\theta }_n(\mathbf {X}) - \theta_0
\bigr),Z \bigr) \nonumber\\
&&\quad\leq\frac{1}{\sqrt{n}} \biggl(2 + \frac
{1}{[i(\theta_0)]^{{3}/{2}}} \biggl[{
\mathrm{E}}\biggl\llvert \frac
{\mathrm
{d}}{\mathrm{d}\theta}\operatorname{ log}f(X_1|
\theta_0)\biggr\rrvert ^3 \biggr] \biggr)
\nonumber
\\[-8pt]
\\[-8pt]
\nonumber
&&\qquad{}+ 2\frac{{\mathrm{E}} (\hat{\theta}_n(\mathbf{X}) -
\theta_0
)^2}{\varepsilon^2} + \frac{1}{\sqrt{n i(\theta_0)}} \biggl\{{\mathrm{E}} \bigl(\bigl
\llvert R_2(\theta _0;\mathbf{X})\bigr\rrvert
|\bigl |\hat{\theta}_n(\mathbf {X}) - \theta_0\bigr|
\leq\varepsilon \bigr)
\\
& &\qquad{}+  \frac{1}{2} \Bigl[{\mathrm{E}} \Bigl( \Bigl(\sup
_{\theta
:|\theta-\theta_0|\leq\varepsilon}\bigl\llvert l^{(3)}(\theta ;\mathbf {X})
\bigr\rrvert \Bigr)^2\big| \bigl|\hat{\theta}_n(
\mathbf{X}) - \theta _0\bigr|\leq\varepsilon \Bigr)
\Bigr]^{{1}/{2}} \bigl[{\mathrm {E}} \bigl(\hat {\theta}_n(
\mathbf{X}) - \theta_0 \bigr)^4 \bigr]^{{1}/{2}}
\biggr\},\nonumber
\end{eqnarray}
where
%
\begin{equation}
\label{R2eq} R_2(\theta_0,\mathbf{x}) = \bigl(
\hat{\theta}_n(\mathbf{x}) - \theta _0\bigr)
\bigl(l''(\theta_0;\mathbf{x}) + n i(
\theta_0) \bigr).
\end{equation}
\end{thmm}
The following lemma is useful for the conditional
expectations in \eqref{boundTHEOREM}; the proof is in the \hyperref
[app]{Appendix}.

%
\begin{Lemma}
\label{Lemmaincreasing}
Let $M \geq0$ be a random variable and $\varepsilon> 0$. For every
continuous function $f$ such that $f(m)$ is increasing and $f(m) \geq
0$, for $m > 0$,
\[
{\mathrm{E}}\bigl[f(M) | M \leq\varepsilon\bigr] \leq{\mathrm{E}}
\bigl[f(M)\bigr].
\]
\end{Lemma}

\begin{pf*}{Proof of Theorem~\ref{Theoremnoncan}}
For the sake of presentation, we drop the subscript $\theta_0$ from the
expectation. The regularity conditions ensure that $0=l'(\hat{\theta
}_n(\mathbf{x});\mathbf{x})$. A second order Taylor expansion
of $l'(\hat{\theta}_n(\mathbf{x});\mathbf{x})$ about
$\theta_0$ gives
%
\begin{equation}
\label{second_order_Taylor} l''(\theta_0;
\mathbf{x}) \bigl(\hat{\theta}_n(\mathbf{x}) -
\theta_0 \bigr) = -l'(\theta_0;
\mathbf{x}) - R_1(\theta _0;\mathbf{x}),
\end{equation}
where
\[
R_1(\theta_0;
{x}) =
\tfrac{1}{2} \bigl(\hat {\theta }_n(\mathbf{x}) -
\theta_0 \bigr)^2 l^{(3)}\bigl(\theta
^{*};\mathbf{x}\bigr)
\]
is the remainder term with $\theta^{*}$ lying between $\hat{\theta
}_n(\mathbf{x})$ and $\theta_0$. The result in \eqref
{second_order_Taylor} gives
\[
-n i(\theta_0) \bigl(\hat{\theta}_n(
\mathbf{x}) - \theta _0 \bigr) = -l'(
\theta_0;\mathbf{x})-R_1(\theta_0;
\mathbf {x}) - \bigl(\hat{\theta}_n(\mathbf{x}) -
\theta_0 \bigr) \bigl[l''(\theta
_0;\mathbf{x}) + n i(\theta_0) \bigr].
\]
As $i(\theta_0)\neq0$
\[
\hat{\theta}_n(\mathbf{x}) - \theta_0 =
\frac
{l'(\theta
_0;\mathbf{x})+R_1(\theta_0;\mathbf{x}) + R_2(\theta
_0,\mathbf{x})}{n i(\theta_0)},
\]
with $R_2(\theta_0,\mathbf{x})$ as in \eqref{R2eq}. For $Z \sim
\mathrm{N}(0,1)$ and $h \in H$ given in \eqref{classfunctions}, we obtain
%
\begin{eqnarray}
&&\bigl\llvert {\mathrm{E}} \bigl[ h \bigl(\bigl(\hat{\theta
}_n(\mathbf{X}) - \theta_0\bigr)\sqrt{n i(
\theta_0)} \bigr) \bigr] - {\mathrm {E}}\bigl[h(Z)\bigr]\bigr\rrvert
\nonumber\\
\label{finalboundnoncanonical2}
&&\quad\leq\biggl\llvert {\mathrm{E}} \biggl[h \biggl(\frac
{l'(\theta_0;\mathbf{X}) + R_1(\theta_0;\mathbf{X}) +
R_2(\theta
_0;\mathbf{X})}{\sqrt{n i(\theta_0)}}
\biggr) - h \biggl(\frac
{l'(\theta_0;\mathbf{X})}{\sqrt{n i(\theta_0)}} \biggr) \biggr]\biggr\rrvert
\\
\label{finalboundnoncanonical1} &&\qquad{} + \biggl\llvert {\mathrm{E}} \biggl[h \biggl(\frac
{l'(\theta_0;\mathbf{X})}{\sqrt{n i(\theta_0)}}
\biggr) \biggr] - {\mathrm{E}}\bigl[h(Z)\bigr]\biggr\rrvert .
\end{eqnarray}
The upper bound for \eqref{finalboundnoncanonical1} is given in
Proposition~\ref{propscore}. To bound \eqref{finalboundnoncanonical2},
note that the term $R_1(\theta_0;\mathbf{X})$ is in general not
uniformly bounded. For ease of presentation, let
\[
C_1= C_1(h,\theta_0;\mathbf{X}) = h
\biggl(\frac{l'(\theta
_0;\mathbf{X}) + R_1(\theta_0;\mathbf{X}) + R_2(\theta
_0;\mathbf{X})}{\sqrt{n i(\theta_0)}} \biggr) - h \biggl(\frac
{l'(\theta_0;\mathbf{X})}{\sqrt{n i(\theta_0)}} \biggr).
\]
For all $x$ the rather crude bound $|C_1|\leq2\|h\|$ is valid. If
$\llvert \hat{\theta}_n(\mathbf{X}) - \theta_0\rrvert \leq
\varepsilon$
then a better bound is available. Hence, we condition on whether $|\hat
{\theta}_n(\mathbf{X}) - \theta_0| > \varepsilon$ or $|\hat
{\theta
}_n(\mathbf{X}) - \theta_0| \leq\varepsilon$, with $\varepsilon
> 0$
such that $(\theta_0 - \varepsilon, \theta_0 + \varepsilon) \subset
\Theta$,
as condition (R3) requires. Moreover, by Markov's inequality
%
\begin{equation}
\label{Chebyshev} \Prob_{\theta_0}\bigl(\bigl|\hat{\theta}_n(
\mathbf{X}) - \theta _0\bigr|>\varepsilon\bigr) \leq
\frac{{\mathrm{E}}[\hat{\theta}_n(\mathbf{X}) - \theta
_0]^2}{\varepsilon^2}.
\end{equation}
Using the law of total expectation,
\begin{eqnarray*}
\bigl|{\mathrm{E}}[C_1]\bigr| &\leq&{\mathrm{E}} \bigl(|C_1|
| \bigl|\hat {\theta }_n(\mathbf{X}) - \theta_0\bigr| >
\varepsilon \bigr)\Prob \bigl(\bigl|\hat {\theta}_n(\mathbf{X}) -
\theta_0\bigr|>\varepsilon \bigr)
\\
&&{}+ {\mathrm{E}} \bigl(|C_1| | \bigl|\hat{\theta}_n(
\mathbf {X}) - \theta _0\bigr|\leq\varepsilon \bigr)\Prob \bigl(\bigl|\hat{
\theta}_n(\mathbf {X})-\theta _0\bigr|\leq\varepsilon
\bigr).
\end{eqnarray*}
Using \eqref{Chebyshev} for\vspace*{-2pt} the first term and a first order Taylor
expansion of $h (\frac{l'(\theta_0;\mathbf{X}) +
R_1(\theta_0;\mathbf{X}) + R_2(\theta_0;\mathbf{X})}{\sqrt{n
i(\theta_0)}}  )$ about $\frac{l'(\theta_0)}{\sqrt{n i(\theta
_0)}}$ for the second term gives
\begin{eqnarray*}
\bigl|{\mathrm{E}}[C_1]\bigr| &\leq&2\|h\|\frac{{\mathrm{E}}
(\hat{\theta
}_n(\mathbf{X}) - \theta_0 )^2}{\varepsilon^2}
\\
&&{} + \biggl\llvert {\mathrm{E}} \biggl(\frac{R_1(\theta
_0,\mathbf{X})
+ R_2(\theta_0,\mathbf{X})}{\sqrt{n i(\theta
_0)}}h'
\bigl(t(\mathbf {X})\bigr)\Big| \bigl|\hat{\theta}_n(
\mathbf{X}) - \theta_0\bigr|\leq \varepsilon \biggr)\biggr\rrvert
\\
&\leq&2\|h\|\frac{{\mathrm{E}} (\hat{\theta
}_n(\mathbf{X})
- \theta_0 )^2}{\varepsilon^2} + \frac{\|h'\|}{\sqrt{n
i(\theta
_0)}}{\mathrm{E}} \bigl(
\bigl\llvert R_2(\theta_0;\mathbf{X})\bigr\rrvert
| \bigl|\hat{\theta}_n(\mathbf{X}) - \theta_0\bigr|
\leq\varepsilon \bigr)
\\
&&{} + \frac{\|h'\|}{\sqrt{n i(\theta_0)}}{\mathrm{E}} \biggl(\frac
{1}{2}\bigl(\hat{
\theta}_n(\mathbf{X}) - \theta_0\bigr)^2
\bigl\llvert l^{(3)}\bigl(\theta ^*;\mathbf{X}\bigr)\bigr\rrvert
\Big| \bigl|\hat{\theta}_n(\mathbf {X}) - \theta_0\bigr|
\leq\varepsilon \biggr),
\end{eqnarray*}
where $t(\mathbf{X})$ lies between $\frac{l'(\theta
_0;\mathbf
{X})}{\sqrt{n i(\theta_0)}}$ and $\frac{l'(\theta_0;\mathbf
{X}) +
R_1(\theta_0;\mathbf{x}) + R_2(\theta_0;\mathbf{X})}{\sqrt{n
i(\theta_0)}}$. Since for $|\hat{\theta}_n(\mathbf{X}) -
\theta_0| \leq\varepsilon$, $|R_1(\theta_0;\mathbf{x})| \leq
\frac
{1}{2}(\hat{\theta}_n(\mathbf{X}) - \theta_0)^2
{\sup}_{\theta:|\theta-\theta_0|\leq\varepsilon}\llvert l^{(3)}(\theta;\mathbf
{X})\rrvert $,
\begin{eqnarray*}
\bigl|{\mathrm{E}}[C_1]\bigr| &\leq&2\|h\|\frac{{\mathrm{E}}
(\hat{\theta
}_n(\mathbf{X}) - \theta_0 )^2}{\varepsilon^2} +
\frac{\|
h'\|
}{\sqrt{n i(\theta_0)}}{\mathrm{E}} \bigl(\bigl\llvert R_2(\theta
_0;\mathbf {X})\bigr\rrvert |\bigl |\hat{\theta}_n(
\mathbf{X}) - \theta_0\bigr| \leq \varepsilon \bigr)
\\
&&{} +\frac{\|h'\|}{2\sqrt{n i(\theta_0)}}{\mathrm{E}} \Bigl[\sup_{\theta:|\theta-\theta_0|\leq\varepsilon}
\bigl\llvert l^{(3)}(\theta ;\mathbf {X})\bigr\rrvert \bigl(\hat{
\theta}_n(\mathbf{X}) - \theta_0
\bigr)^2\big| \bigl|\hat{\theta}_n(\mathbf{X}) -
\theta_0\bigr|\leq \varepsilon \Bigr].
\end{eqnarray*}
The next step is based on the Cauchy--Schwarz inequality and the fact that
%
\begin{equation}
\label{Cauchy} {\mathrm{E}} \bigl[ \bigl(\hat{\theta}_n(
\mathbf{X}) - \theta _0 \bigr)^4| \bigl|\hat{
\theta}_n(\mathbf{X}) - \theta_0\bigr|\leq \varepsilon
\bigr] \leq{\mathrm{E}} \bigl[ \bigl(\hat{\theta}_n(\mathbf {X})
- \theta _0 \bigr)^4 \bigr],
\end{equation}
due to Lemma~\ref{Lemmaincreasing}, giving
%
\begin{eqnarray}
\label{B2noncanonicalbound} \bigl\llvert {\mathrm{E}} [C_1 ]\bigr\rrvert &\leq&2
\|h\|\frac
{{\mathrm{E}} (\hat
{\theta}_n(\mathbf{X}) - \theta_0 )^2}{\varepsilon^2}+ \frac{\|
h'\|}{\sqrt{n i(\theta_0)}} \biggl\{{\mathrm{E}} \bigl(\bigl
\llvert R_2(\theta_0;\mathbf {X})\bigr\rrvert
| \bigl|\hat{\theta}_n(\mathbf{X}) - \theta_0\bigr|
\leq \varepsilon \bigr)
\nonumber
\\
&&{} +  \frac{1}{2} \Bigl[{\mathrm{E}} \Bigl( \Bigl(\sup
_{\theta
:|\theta-\theta_0|\leq\varepsilon}\bigl\llvert l^{(3)}(\theta ;\mathbf {X})
\bigr\rrvert \Bigr)^2\big| \bigl|\hat{\theta}_n(
\mathbf{X}) - \theta _0\bigr|\leq\varepsilon \Bigr)
\Bigr]^{{1}/{2}}
\\
\nonumber
&&{}\times\bigl[{\mathrm {E}} \bigl(\hat {\theta}_n(
\mathbf{X}) - \theta_0 \bigr)^4 \bigr]^{{1}/{2}}
\biggr\}.\qquad
\end{eqnarray}
The result of the theorem is obtained using \eqref{equationpropA1} and
\eqref{B2noncanonicalbound} and the fact that $\|h\| \leq1$ and \mbox{$\|
h'\| \leq1$}.
\end{pf*}

%
\begin{rem}
(1) If $l''(\theta_0;\mathbf{x})\equiv-n i(\theta_0)$
then in \eqref{R2eq}, $R_2(\theta_0;\mathbf{x}) \equiv0$ and the
bound given in Theorem~\ref{Theoremnoncan} simplifies.

(2) The rate of convergence of the Mean Squared Error,
${\mathrm{E}}(\hat{\theta}_n(\mathbf{X}) - \theta_0)^2$, is
$\mathcal{O}
(\frac{1}{n} )$. This result is obtained using that
%
\begin{equation}
\label{MSE} {\mathrm{E}}\bigl(\hat{\theta}_n(\mathbf{X}) -
\theta_0\bigr)^2 = \operatorname{Var}\bigl[\hat {
\theta}_n(\mathbf{X})\bigr] + \operatorname{ bias}^2
\bigl[\hat{\theta }_n(\mathbf{X})\bigr].
\end{equation}
Under the standard asymptotics (from the regularity conditions
(R1)--(R4)) the MLE is asymptotically efficient,
\[
n\operatorname{Var}\bigl[\hat{\theta}_n(\mathbf{X})
\bigr] \mathrel{\mathop{\longrightarrow}_{n \to\infty}} \bigl[i(\theta_0)
\bigr]^{-1},
\]
and hence the variance of the MLE is of order $\frac{1}{n}$. In
addition, from Theorem~\ref{MLEt} the bias of the MLE is of order
$\frac
{1}{\sqrt{n}}$; see also Cox and Snell \cite{Cox}, where no explicit
conditions are
given. Combining these two results and using \eqref{MSE} shows that the
Mean Squared Error of the MLE is of order $\frac{1}{n}$. In the
examples that follow, the remaining terms in the bound are of
order at most $\frac{1}{\sqrt{n}}$.

(3) When the calculation of ${\mathrm{E}} (\llvert \frac
{\mathrm
{d}}{\mathrm{d}\theta}\operatorname{ log}f(X_1|\theta_0)\rrvert ^3 )$ is
awkward, H\"older's inequality can be used, giving ${\mathrm{E}}
(\llvert \frac{\mathrm{d}}{\mathrm{d}\theta}\operatorname{
log}f(X_1|\theta_0)\rrvert ^3 ) \leq [{\mathrm{E}} (\frac{\mathrm{d}}{\mathrm
{d}\theta
}\operatorname{ log}f(X_1|\theta_0) )^4 ]^{{3}/{4}}$.
\end{rem}
%
\section{One-parameter exponential families}
\label{sec:one_parameter}
This section specifies Theorem~\ref{Theoremnoncan} for the distribution
of the MLE for one-parameter exponential family distributions. Many
popular distributions which have the same underlying structure based on
relatively simple properties are exponential families, such as the
normal, Gamma and Laplace distributions. The case of the Poisson
distribution with $\theta\in[0,\infty)$ is treated in Section~\ref
{sec:Poisson}. Generalisations of exponential families can be found in
Lauritzen \cite{Steffen} and Berk \cite{Berk}. The density or
frequency function is of
the form
\[
f(x|\theta) = \operatorname{ exp} \bigl\{ k(\theta)T(x) - A(\theta) +
S(x) \bigr\}\mathbbm{1}_{\{x \in B\}},
\]
where the set $B =  \{ x:f(x|\theta)>0  \}$ is the
support of $X$ and does not depend on $\theta$; $k(\theta)$ and
$A(\theta)$ are functions of the parameter; $T(x)$ and $S(x)$ are
functions only of the data. The choice of the functions $k(\theta)$ and
$T(X)$ is not unique. The case $k(\theta) = \theta$ is the so-called
\emph{canonical case}. In this case, $\theta$ and $T(X)$ are called the
\textit{natural} \textit{parameter} and \textit{natural} \textit
{observation} (Casella and Berger \cite{Casella}).
We make the following assumptions, where (Ass.Ex.1)--(Ass.Ex.3) are
necessary for the existence and uniqueness of the MLE and (A1)--(A4)
follow from the regularity conditions in Section~\ref{sec:intro}.
\begin{enumerate}[(Ass.Ex.1)]
\item[(Ass.Ex.1)]
$\Theta\subset\mathbb{R}$ is open and connected;
\item[(Ass.Ex.2)]
$\lim_{\theta\to\partial\Theta} k(\theta)\sum_{i=1}^{n}T(x_i) - n A(\theta) + \sum_{i=1}^{n}S(x_i)= -\infty$;
\item[(Ass.Ex.3)] We have $k''(\theta)\sum_{i=1}^{n}T(x_i)-n
A''(\theta
) < 0$ at every point
$\theta \in \Theta$ for which it holds that
$k'(\theta) \sum_{i=1}^{n}T(x_i) - n A'(\theta) = 0$;
\item[(A1)]
$k'(\theta) \neq0, \forall\theta\in\Theta$ and $D(\theta) =
\frac
{A'(\theta)}{k'(\theta)}$ is invertible;
\item[(A2)]
$l(\theta;x)$ is thrice continuously differentiable with respect to
$\theta$, meaning that both $k^{(3)}(\theta)$ and $A^{(3)}(\theta)$
exist and they are continuous. In addition, integration of the density
function over $x$ and differentiation with respect to $\theta$ are
three times interchangeable;
\item[(A3)]
for any $\theta_0 \in\Theta$, there exists a positive number
$\varepsilon
$ and a function $M(x)$ (both of which may depend on $\theta_0$) such that
\[
\bigl\llvert k^{(3)}(\theta)T(x) - A^{(3)}(\theta)
\bigr\rrvert \leq M(x)\qquad \forall x \in B, \theta_0 - \varepsilon<
\theta< \theta_0 + \varepsilon,
\]
with ${\mathrm{E}}[M(X)] < \infty$;
\item[(A4)]
$\operatorname{Var}[T(X)] > 0$;
\item[(A5)]
${\mathrm{E}}\llvert T(X) - D(\theta_0)\rrvert ^3$ exists. This
assumption is
required for meaningful bounds.
\end{enumerate}
%

%
\begin{Corollary}
\label{Theoremnoncanexp}
Let $X_1, X_2, \ldots, X_n$ be i.i.d. random variables with the density
or frequency function of a single-parameter exponential family. Assume
that \textup{(A1)--(A5)} are satisfied and that \textup{(Ass.Ex.1)--(Ass.Ex.3)} also hold.
With $Z \sim{\mathrm{N}}(0,1)$, $h \in H$, $R_2(\theta_0;\mathbf{X})$
as in \eqref{R2eq} and also $0 < \varepsilon= \varepsilon(\theta
_0)$ such that $(\theta
_0 - \varepsilon, \theta_0 + \varepsilon)\subset\Theta$ as in
\textup{(A3)}, it holds that
%
\begin{eqnarray}
\label{noncanexponential}
&& d_{bW} \bigl(\sqrt{n i(\theta_0)}\bigl(
\hat{\theta}_n(\mathbf{X}) - \theta_0\bigr),Z \bigr)\nonumber
\nonumber\\
&&\quad\leq\frac{1}{\sqrt{n}} \biggl(2 + \frac
{{\mathrm{E}}|T(X_1)-D(\theta_0)|^3}{ [\operatorname
{Var}[T(X_1)] ]^{{3}/{2}}} \biggr)\nonumber
\\[-8pt]
\\[-8pt]
&&\qquad{}+ 2\frac{{\mathrm{E}} (\hat{\theta}_n(\mathbf
{X}) - \theta
_0 )^2}{\varepsilon^2} + \frac{1}{\sqrt{n i(\theta_0)}} \biggl\{ {\mathrm{E}}
\bigl(\bigl\llvert R_2(\theta_0;\mathbf{X})\bigr
\rrvert |\bigl |\hat{\theta }_n(\mathbf{X}) -
\theta_0\bigr| \leq\varepsilon \bigr)\nonumber
\\
\nonumber
&&\qquad{} +  \frac{1}{2} \Bigl[{\mathrm{E}} \Bigl( \Bigl(\sup
_{\theta
:|\theta-\theta_0|\leq\varepsilon}\bigl\llvert l^{(3)}(\theta ;\mathbf {X})
\bigr\rrvert \Bigr)^2\big| \bigl|\hat{\theta}_n(
\mathbf{X}) - \theta _0\bigr|\leq\varepsilon \Bigr)
\Bigr]^{{1}/{2}} \bigl[{\mathrm {E}} \bigl(\hat {\theta}_n(
\mathbf{X}) - \theta_0 \bigr)^4 \bigr]^{{1}/{2}}
\biggr\}.
\end{eqnarray}
\end{Corollary}
\begin{pf} For the first term of the bound, let
\[
Y_i = Y_i(X_i;\theta
_0) = \biggl(\frac{\mathrm{d}}{\mathrm{d}\theta}\operatorname{
log}f(X_i|\theta
_0)\biggr)\Big/\sqrt{i(\theta_0)},\qquad i=1,2,\ldots,n.
\]
Using Proposition~\ref
{propscore}, we calculate ${\mathrm{E}}\llvert Y_1\rrvert ^3$. Now
\[
\frac{\mathrm{d}}{\mathrm{d}\theta} \log f(X_i|\theta) \Big|_{\theta
=\theta_0} =
k'(\theta_0)T(X_i)-A'(
\theta_0)
\]
yields
\begin{eqnarray*}
{\mathrm{E}}\biggl\llvert \frac{\mathrm{d}}{\mathrm{d}\theta
}\log f(X_i|
\theta_0)\biggr\rrvert ^3 &=& {\mathrm{E}}\bigl\llvert
k'(\theta _0)T(X_i)-A'(\theta
_0)\bigr\rrvert ^3 ={\mathrm{E}}\bigl\llvert
k'(\theta_0) \bigl(T(X_i) - D(\theta
_0)\bigr)\bigr\rrvert ^3
\\
& =& \bigl|k'(\theta_0)\bigr|^3{
\mathrm{E}}\bigl|T(X_i)-D(\theta_0)\bigr|^3\qquad \forall i
\in \{1, 2, \ldots, n \}.
\end{eqnarray*}
In addition, $i(\theta_0) = \operatorname{Var} [\frac{\mathrm
{d}}{\mathrm
{d}\theta} \log f(X_i|\theta_0) ] = [k'(\theta
_0)]^2\operatorname{Var}[T(X_i)] > 0$ from (A1) and (A4). These
quantities can now be
applied to get the first term of the bound in \eqref{noncanexponential}
while the rest of the terms are as in Theorem~\ref{Theoremnoncan}.
\end{pf}
%
\begin{rem}
In the canonical case, $n i(\theta_0) \equiv n A''(\theta_0) \equiv
-l''(\theta_0;\mathbf{x})$. So $R_2(\theta_0;\mathbf{x})
\equiv0$.
\end{rem}
%
\subsection{Example: The exponentially distributed random variable}
\label{sec:examples}
In this section, we consider two examples using the exponential
distribution, first, its canonical form, and then under a change of
parameterisation.
\subsubsection{The canonical case}
\label{subsec:exponentialcan}
In the case of $X_1, X_2, \ldots, X_n$ exponentially distributed,
$\operatorname{Exp}(\theta)$, i.i.d. random variables where $\theta>0$ the probability
density function is
\begin{eqnarray*}
f(x|\theta) = \theta\operatorname{ exp}\{-\theta x\} = \operatorname{
exp}\{\log {\theta} - \theta x\} = \operatorname{ exp} \bigl\{ k(\theta)T(x) - A(
\theta) + S(x) \bigr\}\mathbbm{1}_{\{ {x \in B} \} },
\end{eqnarray*}
where $B = (0,\infty)$, $\theta\in\Theta= (0,\infty)$, $T(x) = -x$,
$k(\theta)=\theta$, $A(\theta) = -\log{\theta}$ and $S(x)=0$. Hence,
$\operatorname{Exp}(\theta)$ is a single-parameter canonical exponential family. Moreover,
\[
l'(\theta;\mathbf{x}) = \frac{n}{\theta} - \sum
_{i=1}^{n}x_i,\qquad
l''(\theta;\mathbf{x}) = -\frac{n}{\theta^2}.
\]
Thus, it is easy to see that the MLE exists, it is unique, equal to
$\hat{\theta}_n(\mathbf{X}) = \frac{1}{\bar{X}}$ and (A1)--(A5) are
satisfied. Corollary~\ref{Theoremnoncanexp} gives
%
\begin{eqnarray}
\label{boundexponential3} d_{bW} \bigl(\sqrt{n i(\theta_0)}\bigl(
\hat {\theta }_n(\mathbf{X}) - \theta_0\bigr),Z
\bigr)&\leq&\frac{4.41456}{\sqrt{n}} + \frac{8(n+2)}{(n-1)(n-2)}
\nonumber
\\[-8pt]
\\[-8pt]
\nonumber
&&{}+ \frac{8\sqrt{n}(n+2)}{(n-1)(n-2)}.
\end{eqnarray}
For $\varepsilon>0$, since $\Theta= (0,\infty)$ simple calculations yield
that $0 < \varepsilon< \theta_0$ to apply (A3) and moreover ${\sup}_{\theta:|\theta
-\theta_0|\leq\varepsilon}\llvert l^{(3)}(\theta; \mathbf
{x})\rrvert  = \frac{2n}{(\theta_0-\varepsilon)^3}$. Choosing $\varepsilon=
\frac{\theta
_0}{2}$, gives that ${\sup}_{\theta:|\theta-\theta_0|\leq
\varepsilon
}\llvert l^{(3)}(\theta;\mathbf{x})\rrvert  = \frac{16n}{\theta
_0^3}$. In addition, since $X_i \sim\operatorname{ Exp}(\theta),
\forall i
\in \{1,2,\ldots, n  \}$ then $\bar{X} \sim\mathrm{
G}(n,n\theta)$, with $G(\alpha, \beta)$ being the Gamma distribution
with shape parameter $\alpha$ and rate parameter $\beta$. Basic\vspace*{1pt}
calculations of integrals show that ${\mathrm{E}}|T(X)-D(\theta_0)|^3
= {\mathrm{E}}\llvert \frac{1}{\theta_0} - X\rrvert ^3 \leq\frac
{2.41456}{\theta
_0^3}$ and
\[
{\mathrm{E}}\bigl[\bigl(\hat{\theta}_n(\mathbf{X}) -
\theta _0\bigr)^2\bigr] = \frac
{(n\theta_0)^2}{(n-1)(n-2)} -
\frac{2n\theta_{0}^{2}}{n-1} + \theta _{0}^{2} = \frac{(n+2)\theta_{0}^{2}}{(n-1)(n-2)}.
\]
Since ${\sup}_{\theta:|\theta-\theta_0|\leq\varepsilon}\llvert l^{(3)}(\theta)\rrvert $ does not depend on the sample, it is not
necessary to use \eqref{Cauchy}. Thus, $\varepsilon= \frac{\theta_0}{2}$
yields the result in \eqref{boundexponential3}.
%
\begin{rem}
(1) The rate of convergence of the bound is $\mathcal
{O}
(\frac{1}{\sqrt{n}} )$. Note also that the bound does not depend
on the value of $\theta_0$.

(2) Note that the\vspace*{-2pt} calculation of ${\mathrm{E}}\llvert \frac
{1}{\theta
_0} - X\rrvert ^3$ requires a significant amount of steps. Therefore,
one could use H\"older's inequality\vspace*{-1pt} with ${\mathrm{E}}\llvert \frac
{1}{\theta
_0} - X\rrvert ^3 \leq [{\mathrm{E}} (\frac{1}{\theta_0}
- X
)^4 ]^{{3}/{4}} = \frac{9^{{3}/{4}}}{\theta_0^3}$ using
the results in pages~70--73 of Kendall and Stuart \cite{DistributionTheory}.
\end{rem}
%
\subsubsection{The non-canonical case}
\label{sec:examplenoncanonical}
Let $X_1, X_2, \ldots, X_n$ be i.i.d. random variables from $\operatorname{Exp}
(\frac{1}{\theta} )$, with p.d.f.
%
\begin{eqnarray}
\label{densityexponentialnon} f(x|\theta) &=& \frac{1}{\theta}\operatorname{ exp} \biggl\{-
\frac
{1}{\theta
}x \biggr\}= \operatorname{ exp} \biggl\{-\operatorname{ log}
\theta- \frac
{1}{\theta}x \biggr\}
\nonumber
\\[-8pt]
\\[-8pt]
\nonumber
&= &\operatorname{ exp} \bigl\{ k(\theta)T(x) - A(\theta) + S(x) \bigr\}
\mathbbm{1}_{\{ {x \in B} \} },
\end{eqnarray}
where $B = (0,\infty)$, $\theta\in\Theta= (0,\infty)$, $T(x) = -x$,
$k(\theta) = \frac{1}{\theta}$, $A(\theta) = \operatorname{
log}\theta$ and $S(x)
= 0$. Again, it is easy to show that the MLE exists, it is unique,
equal to $\hat{\theta}_n(\mathbf{X}) = \bar{X}$ and (A1)--(A5) are
satisfied. For $\varepsilon$ as before and $h \in H$, Corollary~\ref
{Theoremnoncanexp} gives
%
\begin{eqnarray}
\label{boundnoncanexponential3} d_{bW} \bigl(\sqrt{n i(\theta _0)}\bigl(
\hat {\theta}_n(\mathbf{X}) - \theta_0\bigr),Z
\bigr)&\leq&\frac
{4.41456}{\sqrt{n}} + \frac{8}{n}+ \frac{2}{\sqrt{n}}
\nonumber
\\[-8pt]
\\[-8pt]
\nonumber
&&{} + \frac{1}{\sqrt{n}} \biggl(80 \biggl[3 \biggl(\frac{2}{n} + 1
\biggr) \biggr]^{{1}/{2}} \biggr).
\end{eqnarray}
The Mean Squared Error is found to be ${\mathrm{E}}(\hat{\theta
}_n(\mathbf{X}) - \theta_0)^2 = {\mathrm{E}} (\bar{X} -
\theta
_0 )^2 = \frac{\theta_0^2}{n}$.
Also \eqref{densityexponentialnon} gives that $l^{(3)}(\theta
;\mathbf{X}) = -\frac{2n}{\theta^3} + \frac{6}{\theta^4} \sum_{i=1}^{n} X_i = \frac{2n}{\theta^4}(3\hat{\theta}_n(\mathbf
{X}) -
\theta)$ and the triangle inequality yields
\begin{eqnarray*}
\sup_{\theta:|\theta-\theta_0|\leq\varepsilon}\bigl\llvert l^{(3)}(\theta;
\mathbf{X})\bigr\rrvert \leq\sup_{\theta:|\theta
-\theta
_0|\leq\varepsilon} \biggl[\biggl\llvert
\frac{6n\hat{\theta
}_n(\mathbf
{X})}{\theta^4}\biggr\rrvert +\biggl\llvert \frac{2n}{\theta^3}\biggr\rrvert
\biggr] = \frac
{2n}{(\theta_0 - \varepsilon)^4} \bigl(3\hat{\theta}_n(\mathbf {X})+
\theta_0 - \varepsilon \bigr).
\end{eqnarray*}
Therefore,
\begin{eqnarray*}
&& \Bigl[{\mathrm{E}} \Bigl( \Bigl(\sup_{\theta:|\theta
-\theta
_0|\leq\varepsilon}\bigl\llvert
l^{(3)}(\theta;\mathbf{X})\bigr\rrvert \Bigr)^2\big|
\bigl|\hat{\theta}_n(\mathbf{X}) - \theta_0\bigr|\leq
\varepsilon \Bigr) \Bigr]^{{1}/{2}}
\\
&&\quad\leq \biggl[{\mathrm{E}} \biggl( \biggl(\frac{2n}{(\theta_0
- \varepsilon
)^4}\bigl(3\hat{
\theta}_n(\mathbf{X}) + \theta_0 - \varepsilon
\bigr) \biggr)^2\Big| \bigl|\hat{\theta}_n(\mathbf{X})
- \theta_0\bigr|\leq \varepsilon \biggr) \biggr]^{{1}/{2}}
\\
&&\quad\leq\frac{2n}{(\theta_0 - \varepsilon)^4} \bigl[{\mathrm {E}} \bigl( \bigl(3\bigl\llvert
\hat{\theta}_n(\mathbf{X}) - \theta_0\bigr\rrvert +
4\theta _0 - \varepsilon \bigr)^2| \bigl|\hat{
\theta}_n(\mathbf{X}) - \theta _0\bigr|\leq\varepsilon
\bigr) \bigr]^{{1}/{2}}
\\
&&\quad\leq\frac{2n}{(\theta_0 - \varepsilon)^4} \bigl[ (2\varepsilon + 4\theta_0
)^2 \bigr]^{{1}/{2}} = \frac{4n(2\theta_0 +
\varepsilon)}{(\theta_0 - \varepsilon)^4}.
\end{eqnarray*}
The quantity $ [{\mathrm{E}} (\hat{\theta}_n(\mathbf
{X}) - \theta
_0 )^4 ]^{{1}/{2}}$ is calculated using the results in
page 73 and the equations~(3.38), page 70 of Kendall and Stuart \cite
{DistributionTheory} along with the fact that $\hat{\theta
}_n(\mathbf{X}) = \bar{X} \sim G (n, \frac{n}{\theta
_0}
)$, yielding that ${\mathrm{E}}(\hat{\theta}_n(\mathbf{X}) -
\theta_0)^4
= \frac{3\theta_0^4}{n^2} (\frac{2}{n} + 1 )$. Therefore,
\begin{eqnarray*}
&& \Bigl[{\mathrm{E}} \Bigl( \Bigl(
\sup_{\theta:|\theta
-\theta
_0|\leq\varepsilon}\bigl\llvert l^{(3)}(\theta;\mathbf{X})\bigr
\rrvert \Bigr)^2\big| \bigl|\hat{\theta}_n(\mathbf{X})
- \theta_0\bigr|\leq \varepsilon \Bigr) \Bigr]^{{1}/{2}} \bigl[{
\mathrm{E}} \bigl(\hat{\theta }_n(\mathbf{X}) -
\theta_0 \bigr)^4 \bigr]^{{1}/{2}}
\\
&&\quad\leq\frac{4n(2\theta_0 + \varepsilon)}{(\theta_0 - \varepsilon
)^4} \biggl[\frac
{3\theta_0^4}{n^2} \biggl(\frac{2}{n} + 1
\biggr) \biggr]^{{1}/{2}} = \frac{4(2\theta_0 + \varepsilon)}{(\theta_0 - \varepsilon)^4} \biggl[3\theta
_0^4 \biggl(\frac{2}{n} + 1 \biggr)
\biggr]^{{1}/{2}}.
\end{eqnarray*}
To find an upper bound for ${\mathrm{E}} (|R_2(\theta
_0;\mathbf
{X})|||\hat{\theta}_n(\mathbf{X}) - \theta_0|\leq
\varepsilon
 )$,
\begin{eqnarray*}
R_2(\theta_0;\mathbf{X}) &=& \bigl(\hat{
\theta }_n(\mathbf{X}) - \theta_0 \bigr) \biggl(
\frac{n}{\theta_0^2} - \frac
{2n\bar{X}}{\theta_0^3} + \frac{n}{\theta_0^2} \biggr)=\bigl(\hat
{\theta}_n(\mathbf{X}) - \theta_0 \bigr) \biggl(
\frac
{2n}{\theta_0^2} - \frac{2n\bar{X}}{\theta_0^3} \biggr)
\\
&= &-\frac{2n(\hat{\theta}_n({\mathbf{X})} - \theta
_0)^2}{\theta_0^3}.
\end{eqnarray*}
Using Lemma~\ref{Lemmaincreasing} for $f(x)= x^2$ gives
\[
{\mathrm{E}} \bigl[ \bigl(\hat{\theta}_n(\mathbf{X}) - \theta
_0 \bigr)^2| \bigl|\hat{\theta}_n(
\mathbf{X}) - \theta_0\bigr|\leq \varepsilon \bigr] \leq{\mathrm{E}}
\bigl[ \bigl(\hat{\theta}_n(\mathbf {X}) - \theta _0
\bigr)^2 \bigr].
\]
Finally,
\begin{eqnarray*}
{\mathrm{E}} \bigl(\bigl|R_2(
\theta_0;\mathbf{x})\bigr| |\bigl|\hat {\theta}_n(
\mathbf{X}) - \theta_0\bigr|\leq\varepsilon \bigr) &=& {\mathrm{E}}
\biggl(\frac{2n}{\theta_0^3}\bigl(\hat{\theta }_n(\mathbf{X}) -
\theta _0\bigr)^2\Big|\bigl|\hat{\theta}_n(
\mathbf{X}) - \theta_0\bigr|\leq \varepsilon \biggr)
\\
& \leq&\frac{2n}{\theta_0^3}{\mathrm{E}} \bigl[ \bigl(\hat {\theta
}_n(\mathbf{X}) - \theta_0 \bigr)^2
\bigr] = \frac{2}{\theta_0}.
\end{eqnarray*}
Applying now the general result of Corollary~\ref{Theoremnoncanexp} for
$\varepsilon= \frac{\theta_0}{2}$ yields the result in \eqref
{boundnoncanexponential3}.
%
\begin{rem}
\label{remarknoncan}
(1) In this case, the speed of convergence related to the
sample size of the above upper bound is $\mathcal{O}{ (\frac
{1}{\sqrt{n}} )}$ and the bound does not depend on $\theta_0$.

(2) Comparing the upper bound in \eqref
{boundnoncanexponential3} with that in \eqref{boundexponential3} for
the canonical case we see that the first term is the same. However, the
rest of the bound is larger in \eqref{boundnoncanexponential3} than in
\eqref{boundexponential3} $\forall n \in\mathbb{N}$.

(3) In the specific occasion of independent, exponentially
distributed random variables with rate parameter $\frac{1}{\theta_0}$,
the MLE exists, it is unique and equal to $\bar{X}$. Define
$W = \frac{\sqrt{n}(\bar{X} - \theta_0)}{\theta_0} = \frac
{1}{\sqrt
{n}}\sum_{i=1}^{n} Y_i$, where $Y_i = \frac{X_i - \theta_0}{\theta_0}$
are independent, zero mean and unit variance random variables. Also,
${\mathrm{E}}(W) = 0$ and $\operatorname{Var}(W) = \frac{1}{n\theta
_0^2}\sum_{i=1}^{n}\operatorname{Var}(X_i) = 1$. Therefore, \eqref
{equationLemma} can be
used to show
%
\begin{equation}
\label{directsteinexponential} d_{bW} \bigl(\sqrt{n i(\theta_0)}\bigl(
\hat{\theta}_n(\mathbf{X}) - \theta_0\bigr),Z \bigr)
\leq\frac{1}{\sqrt{n}} \biggl(2 + \frac
{1}{\theta
_0^3}{\mathrm{E}}|X_1
- \theta_0|^3 \biggr)\leq\frac{4.41456}{\sqrt{n}}.
\end{equation}
The upper bound given in \eqref{directsteinexponential} as a result of
the direct use of Stein's method is smaller than the upper bound given
in \eqref{boundnoncanexponential3} using the general method explained
in Section~\ref{sec:generalc}. However, in order to apply Stein's
method directly, the quantity $(\hat{\theta}_n(\mathbf{x})-
\theta
_0)\sqrt{n i(\theta_0)}$ is assumed to be a sum of independent random
variables. The general method, on the other hand, gives an upper bound
for~\eqref{interest}, whatever the MLE is, as long as the assumptions
expressed in the beginning of the section hold.
\end{rem}
%
\subsubsection{Empirical results}
\label{subsec:empiricalcanonical}
In this subsection, we study the accuracy of our bounds by simulations.
We start by generating 10\,000 trials of $n$ random independent
observations, $x$, from the exponential distribution. The means for the
canonical and the non-canonical case are equal to 1 and 2,
respectively. We evaluate the MLE, $\hat{\theta}_n(\mathbf{X})$, of
the parameter in each trial, which in turn gives a vector of 10\,000
values. We standardise\vspace*{-2pt} these values and we apply to them the function
$h(x) = \frac{1}{x^2+2}$ with $h \in H$ and $\|h\| = 0.5$, $\|h'\| =
\frac{3\sqrt{1.5}}{16}$ to calculate the expressions in \eqref
{equationpropA11} and \eqref{B2noncanonicalbound}. Finally, we compare
$\llvert {\mathrm{E}} [h (\sqrt{n i(\theta_0)} (\hat
{\theta
}_n(\mathbf{X}) - \theta_0 ) ) ] - {\mathrm
{E}}[h(Z)]\rrvert $ with the sum of the right-hand sides of \eqref
{equationpropA11} and \eqref{B2noncanonicalbound}, using the difference
between their values as a measure of the error. The results presented
in the following tables are based on this particular function $h$ while
the bounded Wasserstein metric is a supremum over a broader class of
test functions, given in \eqref{classfunctions}. Here, ${\mathrm
{E}}[h(Z)] =
0.379$ and the results from the simulations are shown in Tables~\ref
{tableresultexponential} and \ref{tableresultexponentialnoncanonical}.
The tables indicate that $\llvert \hat{\mathrm{E}} [ h
((\hat{\theta
}_n(\mathbf{X}) - \theta_0)\sqrt{n i(\theta_0)}  )
] -
{\mathrm{E}}[h(Z)]\rrvert $, the bound and the error, decrease as the sample
size gets larger. All the values in Table~\ref{tableresultexponential}
are smaller than the respective ones in Table~\ref
{tableresultexponentialnoncanonical}, as expected from Remark~\ref
{remarknoncan}. The bounds are not very good for $n=100$. The reason
might be due to the crude upper bound related to the second term of the
bound in \eqref{noncanexponential}. However, when $n \geq1000$ the
bounds are informative. For the non-canonical case the bounds using
directly Lemma~\ref{Gesinetheorem} are, as expected, much better than
those from the general approach. The bounds are conceptual and better
constraints may be possible.
%
\begin{table}[b]
\caption{Results taken by simulations from the $\operatorname{Exp}(1)$ distribution}\label{tableresultexponential}
\begin{tabular*}{\textwidth}{@{\extracolsep{\fill}}llll@{}}
\hline
$n$ & $\llvert \hat{\mathrm{E}} [ h ((\hat{\theta
}_n(\mathbf
{X}) - \theta_0)\sqrt{n i(\theta_0)}  )  ] - {\mathrm
{E}}[h(Z)]\rrvert $ & Upper bound & Error\\
\hline
\phantom{100\,0}10 & 0.007 & 1.955 & 1.948\\
\phantom{10\,0}100 & 0.002 & 0.336 & 0.334\\
\phantom{10\,}1000 & 0.001& 0.094 & 0.093\\
\phantom{0}10\,000 & 0.0002 & 0.029 & 0.0288\\
100\,000 & 0.0001 & 0.009 & 0.0089
\\\hline
\end{tabular*}
\end{table}
%
\begin{table}[t]
\caption{Results taken by simulations from the $\operatorname
{Exp}(0.5)$ distribution
treated as a non-canonical exponential family}\label{tableresultexponentialnoncanonical}
\begin{tabular*}{\textwidth}{@{\extracolsep{\fill}}lld{2.3}ll@{}}
\hline
$n$ & $\llvert \hat{\mathrm{E}} [ h ((\hat{\theta
}_n(\mathbf
{X}) - \theta_0)\sqrt{n i(\theta_0)}  )  ] - {\mathrm
{E}}[h(Z)]\rrvert $ & \multicolumn{1}{c}{Bound} & Error & Bound using Lemma~\ref
{Gesinetheorem}\\
\hline
\phantom{100\,0}10 & 0.004\phantom{0} & 11.888 & 11.884 & 0.321\\
\phantom{10\,0}100 &  0.003\phantom{0} &  3.401 &  \phantom{0}3.398 & 0.101\\
\phantom{10\,}1000 &  0.002\phantom{0} &  1.058 &  \phantom{0}1.056 & 0.032\\
\phantom{0}10\,000 &  0.001\phantom{0} &  0.333 &  \phantom{0}0.332 & 0.010\\
100\,000 & 0.0005 & 0.105 & \phantom{0}0.1045 & 0.003\\
\hline
\end{tabular*}
\end{table}

\section{Discrete distributions: The boundary issue}
\label{sec:perturbationc}
In this section, we use a perturbation method for any discrete
distribution that faces the problem of the MLE having positive
probability of being on the boundary of the parameter space. We also
illustrate the perturbation for the specific example of the Poisson
distribution.
\subsection{The perturbation approach}
A perturbation method based on a perturbation function, should be such
that first of all, the function should perturb the quantity of interest
in a way that ensures it will be interior to its domain. The second
requirement is that the perturbed quantity should be as close as
possible to the initial quantity. Let $X$ be a random variable with
support $B$, the connected closed (semi-closed) interval $[a,b]$
($(a,b]$ or $[a,b)$), where\vspace*{-2pt} $-\infty< a < b < \infty$. For $0 <
\varepsilon< \frac{b-a}{2}$, we are looking for a perturbation function,
$q: B \rightarrow\overset{\circ}{B}$ (where in this case, $\overset
{\circ}{B}$ denotes the interior of the set $B$) with $q(x) = kx+d$,
such that:
\begin{enumerate}[(1)]
\item[(1)] $q(a) = a + \varepsilon$ and $q(b) = b - \varepsilon$.
\item[(2)] ${\sup}_{x} |q(x) - x|$ is minimum, $x \in B$.
\end{enumerate}
Solving this problem for $k$ and $d$, gives $k = 1 - \frac
{2\varepsilon
}{b-a}$ and $d = \varepsilon+ \frac{2a}{b-a}\varepsilon$. There is
only one
solution, which is minimal. Thus, the second requirement is also
satisfied. Choose $\varepsilon= \varepsilon(n) = \frac{c}{n}$ and $0
< c <
\frac{n(b-a)}{2}$. Finally, the perturbation function is
%
\begin{equation}
\label{perturbationfunction} q(x) = x + \frac{c}{n} -\frac{2c}{n} \biggl(
\frac{x-a}{b-a} \biggr),\qquad  x \in B,0 < c < \frac{n(b-a)}{2}.
\end{equation}
In the case where $B = (-\infty,b]$ or $B = [a,\infty)$, then $q(x) =
x-\frac{c}{n}$ or $q(x) = x+\frac{c}{n}$, respectively.

Assuming existence and uniqueness of the MLE, $\hat{\theta
}_n(\mathbf{X})$, for the parameter $\theta_0$, of a discrete
distribution with parameter space as in the previous paragraph, the aim
is to find an upper bound on
\[
d_{bW} \bigl(\sqrt{n}\bigl(\hat{\theta}_n(
\mathbf{X}) - \theta _0\bigr),K \bigr),
\]
where $K \sim{\mathrm{N}} (0,\frac{1}{i(\theta_0)} )$.
Note that
$\mathrm{N}(0,0)$ is point mass at 0. The quantity we will bound is
not exactly
the one shown in \eqref{interest} because the Expected Fisher
Information number might not exist or not be finite when $\theta_0$
lies on the boundary of the parameter space. For this purpose, we will
use the perturbation function in \eqref{perturbationfunction} for both
the parameter and the data.

First, we introduce some notations. For $S$ being the discrete sample
space, let $a := \inf{\Theta}$, $b := \sup{\Theta}$, $S_1 :=
\inf{S}$, $S_p := \sup{S}$ and $0 < c_1 < \frac{n(b-a)}{2}$, $0 <
c_2< \frac{n(S_p-S_1)}{2}$. In addition, $\theta_0^*=\theta_0 +
\frac
{c_1}{n} - \frac{2c_1}{n} (\frac{\theta_0-a}{b-a} )$ is the
perturbed parameter and
%
\begin{equation}
\label{perturbed_DATA} q(x_i)= x_i + \frac{c_2}{n} -
\frac{2c_2}{n} \biggl(\frac{x_i - S_1}{S_p
- S_1} \biggr)
\end{equation}
is the perturbed data. The perturbed MLE is denoted by $\hat{\theta
}_n^*(\mathbf{x}):=\hat{\theta}_n(\mathbf{x})|_{\mathbf{x}
= q(\mathbf{x})}$. Also,
\begin{eqnarray*}
 l'\bigl(\theta_0^*;q(\mathbf{x})\bigr)
&:= &l'(\theta ;\mathbf {x})\Big| \mathop{{}_{\theta= \theta_0^*}}_{ \mathbf{x}= q(\mathbf
{x})},\qquad
l''\bigl(\theta_0^*;q(\mathbf{x})
\bigr) := l''(\theta ;\mathbf{x}) \Big|\mathop{{}_{\theta= \theta_0^*}}_{ \mathbf{x}=
q(\mathbf{x})},
\\
 l^{(3)}\bigl(\theta;q(\mathbf{x})\bigr)& =&
l^{(3)}(\theta ;\mathbf{x}) \big|_{{\mathbf{x}= q(\mathbf{x})}}.
\end{eqnarray*}
For ease of presentation, abbreviate $Y_i = \frac{l'(\theta
_0^*;q(X_i))}{\sqrt{n}i(\theta_0^*)}$, $i \in \{1,\ldots,n
 \}$ while $w_1:=w_1(n,\theta_0^*)$ and $w_2:=w_2(n,\theta
_0^*)$ are its expectation and variance, respectively.

%
\begin{thmm}
\label{generalperbound}
Let $X_1, X_2, \ldots, X_n$ be i.i.d. random variables from a
single-parameter discrete distribution with parameter space the
connected, closed or semi-closed interval $\Theta\subset\mathbb{R}$
and discrete sample space $S$. Assume that $i(\theta_0) > 0$ and let
$\frac{1}{i(\theta_0)} = 0$ to be the continuous extension of $\frac
{1}{i(\theta)}$ to $\theta\rightarrow\theta_0$ when $\theta_0$ is
such that $i(\theta_0)$ does not exist or it is equal to infinity. Let\vspace*{-2pt}
$h \in H$ and $0 <\varepsilon= \varepsilon(\theta_0^*)$ such that
$(\theta^*_0 - \varepsilon, \theta
^*_0 +
\varepsilon)\subset\overset{\circ}{\Theta}$.
Then
%
\begin{eqnarray}
\label{generalperbound1}
&& d_{bW} \bigl(\sqrt{n}\bigl(\hat{
\theta}_n(\mathbf{X}) - \theta _0\bigr),K \bigr)\nonumber\\
&&\quad\leq
\frac{c_1}{\sqrt{n}}\biggl\llvert 1-2 \biggl(\frac{\theta
_0 -
a}{b-a} \biggr)\biggr
\rrvert + \sqrt{n} {\mathrm{E}}\bigl\llvert \hat{\theta }_n(
\mathbf {X}) - \hat{\theta}_n^*(\mathbf{X})\bigr\rrvert
\nonumber\\
&&\qquad{}+ \biggl[\biggl\llvert 1 - \frac{1}{\sqrt{w_2n i(\theta
_0)}}\biggr\rrvert \sqrt{n
w_2 + (n w_1)^2} + \frac{\sqrt{n}|w_1|}{\sqrt{w_2 i(\theta
_0)}}
\biggr]\mathbbm{1} \biggl\{\frac{1}{i(\theta_0)} > 0 \biggr\}
\nonumber
\\[-8pt]
\\[-8pt]
&&\qquad{}+ \frac{1}{\sqrt{n}} \biggl(2 + \frac{1}{(w_2)^{{3}/{2}}}{\mathrm{E}}|Y_1
- w_1|^3 \biggr)\mathbbm{1} \biggl\{\frac{1}{i(\theta_0)}
> 0 \biggr\}+ 2\frac{{\mathrm{E}} (\hat{\theta}^*_n(\mathbf{X}) -
\theta
^*_0 )^2}{\varepsilon^2}\nonumber
\\
\nonumber
&&\qquad{} + \frac{1}{\sqrt{n} i(\theta^*_0)}
\biggl\{{\mathrm {E}} \bigl(\bigl|\bigl(\hat{\theta
}^*_n(\mathbf{x}) - \theta^*_0\bigr)
\bigl[l''\bigl(\theta^*_0;q(\mathbf
{x})\bigr) + n i\bigl(\theta^*_0\bigr)\bigr]\bigr||\bigl|\hat{
\theta}^*_n(\mathbf{X}) - \theta ^*_0\bigr|\leq
\varepsilon \bigr)
\\
\nonumber
&&\qquad{} + \frac{1}{2} \Bigl[{\mathrm{E}} \Bigl( \Bigl(\sup
_{\theta
:|\theta-\theta^*_0|\leq\varepsilon}\bigl\llvert l^{(3)}\bigl(\theta ;q(\mathbf
{X})\bigr)\bigr\rrvert \Bigr)^2\big| \bigl|\hat{\theta}^*_n(
\mathbf{X}) - \theta ^*_0\bigr|\leq\varepsilon \Bigr)
\Bigr]^{{1}/{2}} \bigl[{\mathrm {E}} \bigl(\hat {\theta}^*_n(
\mathbf{X}) - \theta^*_0 \bigr)^4
\bigr]^{{1}/{2}} \biggr\}.
\end{eqnarray}
\end{thmm}
\begin{pf}
\textit{Step} 1: \textit{Perturbation of} $\theta_0$. Using the triangle
inequality and then a first order Taylor expansion of $h(\sqrt{n}(\hat
{\theta}_n(\mathbf{X}) - \theta_0))$ about $\sqrt{n}(\hat
{\theta
}_n(\mathbf{X}) - \theta_0^*)$ gives
%
\begin{eqnarray}
\label{generalperdelta1}
&&\bigl|{\mathrm{E}}\bigl[h\bigl(\sqrt{n}\bigl(\hat{
\theta}_n(\mathbf{X}) - \theta _0\bigr)\bigr)\bigr]
- {\mathrm{E}}\bigl[h(K)\bigr]\bigr|\nonumber\\
&&\quad \leq\bigl|{\mathrm{E}}\bigl[h\bigl(\sqrt{n}\bigl(\hat
{\theta }_n(\mathbf{X}) - \theta^*_0\bigr)\bigr)
\bigr] - {\mathrm{E}}\bigl[h(K)\bigr]\bigr|
\nonumber\\
&&\qquad{}+\bigl |{\mathrm{E}}\bigl[h\bigl(\sqrt{n}\bigl(\hat{\theta}_n(
\mathbf {X}) - \theta _0\bigr)\bigr) - h\bigl(\sqrt{n}\bigl(\hat{
\theta}_n(\mathbf{X}) - \theta^*_0\bigr)\bigr)
\bigr]\bigr|
\\
&&\quad \leq\bigl|{\mathrm{E}}\bigl[h\bigl(\sqrt{n}\bigl(\hat{\theta
}_n(\mathbf{X}) - \theta^*_0\bigr)\bigr)\bigr] - {
\mathrm{E}}\bigl[h(K)\bigr]\bigr| + \sqrt{n}\bigl\|h'\bigr\|{\mathrm {E}}\bigl|
\theta^*_0 - \theta_0\bigr|\nonumber
\\
&&\quad= \bigl|{\mathrm{E}}\bigl[h\bigl(\sqrt{n}\bigl(\hat{\theta}_n(
\mathbf{X}) - \theta^*_0\bigr)\bigr)\bigr] - {\mathrm{E}}
\bigl[h(K)\bigr]\bigr| + \frac{\|h'\|c_1}{\sqrt{n}}\biggl\llvert 1 - 2 \biggl(
\frac{\theta
_0-a}{b-a} \biggr)\biggr\rrvert .\nonumber
\end{eqnarray}
\textit{Step} 2: \textit{Perturbation of the MLE}. To perturb the MLE, we
perturb the data. The perturbed data is denoted by $q(\mathbf{x}) =
(q(x_1), q(x_2), \ldots, q(x_n))$, with $q(x_i)$ given in \eqref
{perturbed_DATA}. This construction ensures that the MLE evaluated at
$q(\mathbf{x})$ is not on the boundary of the parameter space.
Following the same process as in \eqref{generalperdelta1}, using the
triangle inequality and a first order Taylor expansion of $h(\sqrt
{n}(\hat{\theta}_n(\mathbf{X}) - \theta^*_0))$ about $\sqrt
{n}(\hat
{\theta}^*_n(\mathbf{X}) - \theta^*_0)$ gives
%
\begin{eqnarray}
\label{generalperdelta2}
& &\bigl|{\mathrm{E}}\bigl[h\bigl(\sqrt{n}\bigl(\hat{
\theta}_n(\mathbf {X}) - \theta ^*_0\bigr)\bigr)
\bigr] - {\mathrm{E}}\bigl[h(K)\bigr]\bigr|\nonumber\\
&&\quad \leq\bigl|{\mathrm{E}}\bigl[h\bigl(\sqrt{n}
\bigl(\hat {\theta }^*_n(\mathbf{X}) - \theta^*_0
\bigr)\bigr)\bigr] - {\mathrm{E}}\bigl[h(K)\bigr]\bigr|
\nonumber
\\[-8pt]
\\[-8pt]
\nonumber
&&\qquad{}+ \bigl|{\mathrm{E}}\bigl[h\bigl(\sqrt{n}\bigl(\hat{\theta}_n(
\mathbf {X}) - \theta ^*_0\bigr)\bigr) - h\bigl(\sqrt{n}\bigl(
\hat{\theta}^*_n(\mathbf{X}) - \theta ^*_0\bigr)
\bigr)\bigr]\bigr|
\\
&&\quad\leq\bigl|{\mathrm{E}}\bigl[h\bigl(\sqrt{n}\bigl(\hat{\theta}^*_n(
\mathbf{X}) - \theta ^*_0\bigr)\bigr)\bigr] - {\mathrm{E}}
\bigl[h(K)\bigr]\bigr| + \sqrt{n}\bigl\|h'\bigr\|{\mathrm{E}}\bigl|\hat {\theta
}_n(\mathbf{X}) - \hat{\theta}^*_n(
\mathbf{X})\bigr|.\nonumber
\end{eqnarray}
\textit{Step} 3: \textit{The final bound}. It remains to bound
\[
\bigl|{\mathrm{E}}\bigl[h\bigl(\sqrt{n}\bigl(\hat{\theta}^*_n(
\mathbf{X}) - \theta^*_0\bigr)\bigr)\bigr] - {\mathrm{E}}
\bigl[h(K)\bigr]\bigr|.
\]
Since both $\theta_0^*$ and $\hat{\theta}^*_n(\mathbf{x})$ are
interior to $\Theta$, a second-order Taylor expansion of $l'(\hat
{\theta
}^*_n(\mathbf{x});q(\mathbf{x}))$ about $\theta^*_0$ yields
%
\begin{equation}
\label{Taylorgeneralper} 0 = l'\bigl(\theta_0^*;q(
\mathbf{x})\bigr) + \bigl(\hat{\theta }^*_n(\mathbf {x}) -
\theta^*_0\bigr)l''\bigl(
\theta^*_0;q(\mathbf{x})\bigr) + R_1\bigl(\theta
^*_0;q(\mathbf{x})\bigr),
\end{equation}
where, similarly as in Section~\ref{sec:generalc},
\[
R_1\bigl(\theta^*_0;q(\mathbf{x})\bigr) =
\tfrac{1}{2}\bigl(\hat{\theta }^*_n(\mathbf{x}) -
\theta^*_0\bigr)^2l^{(3)}\bigl(\tilde{\theta
};q(\mathbf{x})\bigr)
\]
with
\[
l^{(3)}\bigl(\tilde{\theta};q(\mathbf{x})\bigr) = l^{(3)}(
\theta ;\mathbf {x}) \Big|\mathop{{}_{\theta= \tilde{\theta}}}_{ \mathbf{x}=
q(\mathbf{x})}
\]
for $\tilde{\theta}$ between $\hat{\theta}^*_n(\mathbf{x})$ and
$\theta^*_0$. A simple rearrangement of the terms in \eqref
{Taylorgeneralper}, leads to $\hat{\theta}^*_n(\mathbf{x}) -
\theta
^*_0 = \frac{-l'(\theta^*_0;g(\mathbf{x})) - R_1(\theta
^*_0;g(\mathbf{x}))}{l''(\theta^*_0;g(\mathbf{x}))}$. Since, in
general $l''(\theta_0^*;q(\mathbf{x})) \neq-n i(\theta_0^*)$,
using the results in the proof of Theorem~\ref{Theoremnoncan} gives
\[
\hat{\theta}^*_n(\mathbf{x}) - \theta^*_0
= \frac
{l'(\theta^*_0;q(\mathbf{x})) + R_1(\theta^*_0;q(\mathbf
{x})) +
R_2(\theta^*_0;q(\mathbf{x}))}{n i(\theta^*_0)},
\]
where
\[
R_2\bigl(\theta^*_0;q(\mathbf{x})\bigr) = \bigl(
\hat{\theta}^*_n(\mathbf {x}) - \theta^*_0\bigr)
\bigl[l''\bigl(\theta^*_0;q(\mathbf{x})
\bigr) + n i\bigl(\theta^*_0\bigr)\bigr].
\]
Using that $q(\mathbf{X}) = (q(X_1), q(X_2), \ldots, q(X_{n}))$,
the triangle inequality gives
%
\begin{eqnarray}
\label{A1generalper} &&\bigl|{\mathrm{E}}\bigl[h\bigl(\sqrt{n}\bigl(\hat{
\theta}^*_n(\mathbf{X}) - \theta^*_0\bigr)\bigr)
\bigr] - {\mathrm{E}}\bigl[h(K)\bigr]\bigr|\nonumber\\
&&\quad \leq\biggl\llvert {\mathrm{E}} \biggl[h
\biggl(\frac
{l'(\theta
^*_0;q(\mathbf{X}))}{\sqrt{n} i(\theta^*_0)} \biggr) \biggr] - {\mathrm{E}}\bigl[h(K)\bigr]\biggr
\rrvert
\\
\nonumber
& &\qquad{}+ \biggl\llvert {\mathrm{E}} \biggl[h \biggl(\frac{l'(\theta
^*_0;q(\mathbf
{X})) + R_1(\theta^*_0;q(\mathbf{X})) + R_2(\theta
^*_0;q(\mathbf
{X}))}{\sqrt{n} i(\theta^*_0)}
\biggr) - h \biggl(\frac{l'(\theta
^*_0;q(\mathbf{X}))}{\sqrt{n} i(\theta^*_0)} \biggr) \biggr]\biggr\rrvert .
\end{eqnarray}

(A) To find an upper bound on the first quantity on the
right-hand side of \eqref{A1generalper} using Lemma~\ref
{Gesinetheorem}, note that
\[
\frac{l'(\theta_0^*;q(\mathbf{X}))}{\sqrt{n} i(\theta
_0^*)} = \sum_{i=1}^{n}Y_i,\qquad
\mbox{where } Y_i = \frac{l'(\theta
_0^*;q(X_i))}{\sqrt{n}i(\theta_0^*)}.
\]
Denote by $w_1:=w_1(n)$ and $w_2:=w_2(n)$ the expectation and the
variance of $Y_i, i=1,2,\ldots,n$, respectively. These
quantities depend on the sample size and on the perturbed values
($\theta^*_0$ and $q(x_i)$). Define $\tilde{Y}_i = \frac{Y_i -
w_1}{\sqrt{w_2 i(\theta_0)}}, \forall i \in\{1,2, \ldots, n\}$ with
${\mathrm{E}}(\tilde{Y}_i) = 0$ and $\operatorname{Var}(\tilde
{Y}_i) = \frac
{1}{i(\theta_0)}$. As\vspace*{-1pt} a consequence of $X_1, X_2, \ldots, X_n$ being
i.i.d. random variables, $\tilde{Y}_1, \tilde{Y}_2, \ldots, \tilde
{Y}_n$ are i.i.d. random variables too. Using the triangle inequality
and that
\[
\frac{1}{\sqrt{n}}\sum_{i=1}^{n}
\tilde{Y}_i = \frac{1}{\sqrt{w_2n
i(\theta_0)}} \biggl(\frac{l'(\theta_0^*;q(\mathbf{X}))}{\sqrt{n
i(\theta_0^*)}} - n
w_1 \biggr)
\]
gives
%
\begin{eqnarray}
\label{firsttermgeneralper}
&&\biggl\llvert {\mathrm{E}} \biggl[h \biggl(\frac{l'(\theta^*_0;q(\mathbf
{X}))}{\sqrt
{n} i(\theta^*_0)}
\biggr) \biggr] - {\mathrm{E}}\bigl[h(K)\bigr]\biggr\rrvert
\nonumber\\
&&\quad\leq\biggl\llvert {\mathrm{E}} \biggl[h \biggl(\frac{1}{\sqrt{w_2n
i(\theta
_0)}} \biggl(
\frac{l'(\theta_0^*;q(\mathbf{X}))}{\sqrt{n}
i(\theta
_0^*)} - n w_1 \biggr) \biggr) \biggr] - {\mathrm{E}}
\bigl[h(K)\bigr]\biggr\rrvert
\\
\nonumber
&&\qquad{}+\biggl\llvert {\mathrm{E}} \biggl[h \biggl(\frac{l'(\theta
_0^*;q(\mathbf
{X}))}{\sqrt{n} i(\theta_0^*)}
\biggr) - h \biggl(\frac{1}{\sqrt{w_2n
i(\theta_0)}} \biggl(\frac{l'(\theta_0^*;q(\mathbf{X}))}{\sqrt{n}
i(\theta_0^*)} - n
w_1 \biggr) \biggr) \biggr]\biggr\rrvert .
\end{eqnarray}
The first term of the bound in \eqref{firsttermgeneralper} will be
bounded using Lemma~\ref{Gesinetheorem} with $W = \frac
{1}{\sqrt{n}}\sum_{i=1}^{n} \tilde{Y}_i$. Thus,
%
\begin{eqnarray}
\label{upboundfirsttermgeneralper}
&&\biggl\llvert {\mathrm{E}} \biggl[h \biggl(
\frac{1}{\sqrt{w_2n
i(\theta
_0)}} \biggl(\frac{l'(\theta_0^*;q(\mathbf{X}))}{\sqrt{n}
i(\theta
_0^*)} - n w_1 \biggr)
\biggr) \biggr] - {\mathrm{E}}\bigl[h(K)\bigr]\biggr\rrvert
\nonumber
\\[-8pt]
\\[-8pt]
\nonumber
&&\quad\leq\frac{\|h'\|}{\sqrt{n}} \bigl(2 + \bigl[i(\theta_0)
\bigr]^{{3}/{2}}{\mathrm{E}}\llvert \tilde{Y}_1\rrvert
^3 \bigr) = \frac{\|h'\|}{\sqrt{n}} \biggl(2 + \frac{1}{(w_2)^{{3}/{2}}}{
\mathrm{E}}|Y_1 - w_1|^3 \biggr).
\end{eqnarray}
For the second term of the upper bound in \eqref{firsttermgeneralper} a
first-order Taylor expansion and the Cauchy--Schwarz inequality yield
%
\begin{eqnarray}
\label{upboundsecondtermgeneralper}
&&\biggl\llvert {\mathrm{E}} \biggl[h \biggl(\frac{l'(\theta_0^*;q(\mathbf
{X}))}{\sqrt
{n} i(\theta_0^*)}
\biggr) - h \biggl(\frac{1}{\sqrt{w_2n i(\theta
_0)}} \biggl(\frac{l'(\theta_0^*;q(\mathbf{X}))}{\sqrt{n}
i(\theta
_0^*)} - n
w_1 \biggr) \biggr) \biggr]\biggr\rrvert
\nonumber\\
&&\quad\leq\bigl\|h'\bigr\|\biggl\llvert 1 - \frac{1}{\sqrt{w_2n i(\theta
_0)}}\biggr
\rrvert {\mathrm{E}}\biggl\llvert \frac{l'(\theta_0^*;q(\mathbf{X}))}{\sqrt
{n} i(\theta
_0^*)}\biggr\rrvert +
\frac{\|h'\|\sqrt{n}|w_1|}{\sqrt{w_2 i(\theta
_0)}}
\nonumber
\\
& &\quad\leq\bigl\|h'\bigr\|\biggl\llvert 1 - \frac{1}{\sqrt{w_2n i(\theta
_0)}}\biggr
\rrvert \biggl(\operatorname{Var} \biggl(\frac{l'(\theta_0^*;q(\mathbf
{X}))}{\sqrt{n}
i(\theta_0^*)} \biggr) +
\frac{ [{\mathrm{E}}(l'(\theta
_0^*;q(\mathbf
{X}))) ]^2}{n [i(\theta_0^*)]^2} \biggr)^{{1}/{2}}\\
&&\qquad{} + \frac
{\|
h'\|\sqrt{n}|w_1|}{\sqrt{w_2 i(\theta_0)}}
\nonumber\\
& &\quad=\bigl \|h'\bigr\| \biggl[\biggl\llvert 1 - \frac{1}{\sqrt{w_2n i(\theta
_0)}}
\biggr\rrvert \sqrt{n w_2 + (n w_1)^2} +
\frac{\sqrt{n}|w_1|}{\sqrt
{w_2 i(\theta_0)}} \biggr].\nonumber
\end{eqnarray}
When $\frac{1}{i(\theta_0)} = 0$\vspace*{1pt} then $\tilde{Y}_i = 0,\forall i \in
 \{1, 2, \ldots, n \}$ and by following the above
process, the first term on the right-hand side of \eqref{A1generalper}
is equal to zero.

(B) To complete the proof, it remains to find an upper bound
for the second term on the right-hand side of \eqref{A1generalper}. The
idea is the same as the one used for \eqref{B2noncanonicalbound}. We
condition on whether $|\hat{\theta}^*_n(\mathbf{X}) - \theta
^*_0| >
\varepsilon$ or $|\hat{\theta}^*_n(\mathbf{X}) - \theta^*_0|
\leq
\varepsilon$, where now $\varepsilon= \varepsilon(\theta^*_0)$ and
$0 <\varepsilon
$
$(\theta^*_0 - \varepsilon, \theta
^*_0 +
\varepsilon)\subset\overset{\circ}{\Theta}$.
Following
the same process as in Section~\ref{sec:generalc} yields
%
\begin{eqnarray}
\label{Bgeneralperbound} &&\biggl\llvert {\mathrm{E}} \biggl[h \biggl(\frac{l'(\theta^*_0;g(\mathbf
{X})) +
R_1(\theta^*_0;g(\mathbf{X})) + R_2(\theta^*_0;g(\mathbf
{X}))}{\sqrt{n} i(\theta^*_0)}
\biggr) - h \biggl(\frac{l'(\theta
^*_0;g(\mathbf{X}))}{\sqrt{n} i(\theta^*_0)} \biggr) \biggr]\biggr\rrvert
\nonumber\\
&&\quad\leq2\|h\|\frac{{\mathrm{E}} (\hat{\theta
}^*_n(\mathbf
{X}) - \theta^*_0 )^2}{\varepsilon^2} + \frac{\|h'\|}{\sqrt{n}
i(\theta^*_0)} \biggl\{{
\mathrm{E}} \bigl(\bigl|R_2\bigl(\theta ^*_0;g(\mathbf{X})
\bigr)\bigr||\bigl|\hat {\theta}^*_n(\mathbf{X}) -
\theta^*_0\bigr|\leq\varepsilon \bigr)
\\
\nonumber
& &\qquad{}+ \frac{\|h'\|}{2} \Bigl[{\mathrm{E}} \Bigl( \Bigl(\sup
_{\theta
:\bigl|\theta-\theta^*_0\bigr|\leq\varepsilon}\bigl\llvert l^{(3)}\bigl(\theta ;g(\mathbf
{X})\bigr)\bigr\rrvert \Bigr)^2\big| \bigl|\hat{\theta}^*_n(
\mathbf{X}) - \theta ^*_0\bigr|\leq\varepsilon \Bigr)
\Bigr]^{{1}/{2}} \bigl[{\mathrm {E}} \bigl(\hat {\theta}^*_n(
\mathbf{X}) - \theta^*_0 \bigr)^4
\bigr]^{{1}/{2}} \biggr\}.
\end{eqnarray}
Combining \eqref{generalperdelta1}, \eqref{generalperdelta2}, \eqref
{upboundfirsttermgeneralper}, \eqref{upboundsecondtermgeneralper} and
\eqref{Bgeneralperbound} and the fact that $\|h\|\leq1$, $\|h'\|\leq
1$ gives the result in \eqref{generalperbound1}.
\end{pf}
%
\begin{rem}
\label{remarkgeneral}
(1) In order for the above bound to approach zero as the
sample size, $n$, increases we require that ${\mathrm{E}}\llvert \hat
{\theta
}_n(\mathbf{X}) - \hat{\theta}_n^*(\mathbf{X})\rrvert  =
o
(\frac{1}{\sqrt{n}} )$.

(2) When both endpoints of the parameter space are not finite,
then parameter perturbation is not necessary. In the case where one of
the two endpoints of the now semi-closed parameter space is infinite,
then it suffices to change the form of the perturbed parameter, which
now becomes
\begin{eqnarray*}
\theta_0^*& = &\theta_0 - \frac{c_1}{n}\qquad
\mbox{if the left endpoint is equal to $-\infty$},
\\
\theta_0^*& =& \theta_0 + \frac{c_1}{n}\qquad
\mbox{if the right endpoint is equal to $\infty$}.
\end{eqnarray*}
The same holds regarding the sample space and the relevant perturbation
of the data.
\end{rem}
%
\subsection{Example: The Poisson distribution}
\label{sec:Poisson}
In this subsection, we consider the Poisson distribution with parameter
$\theta\in\Theta= [0,\infty)$. The value $\theta= 0$ must be in the
parameter space in order for the MLE, $\hat{\theta}_n(\mathbf
{X}) =
\bar{X}$, to exist and to be unique. The $\operatorname{Poisson}(\theta)$ distribution
with the aforementioned parameter space is not a single-parameter
exponential family. When $\theta= 0$ is included in the parameter
space the requirements of an exponential family are not satisfied as
the set of values $x$ for which the relevant probability mass function
\[
f(x|\theta) = \frac{\mathrm{ e}^{-\theta}\theta^x}{x!},\qquad \theta \in[0,\infty), x\in
\mathbb{Z}^+_0
\]
is positive, is different for $\theta= 0$ than for any other value of
the parameter $\theta$; the support of the distribution depends on the
parameter. Following the steps of the proof of Theorem~\ref
{generalperbound}, using also H\"{o}lder's inequality for the third
absolute moment in the third term of the bound in \eqref
{generalperbound1} and taking $0 < c = c_1 =c_2$, which minimizes the
bound, gives the next result.

%
\begin{Corollary}
\label{Poissonbound}
Let $X_1, X_2, \ldots, X_n$ be i.i.d. random variables which follow the
$\operatorname{Poisson}(\theta_0)$ distribution, with $\theta_0 \in[0,\infty)$. For
$K \sim{\mathrm{N}} (0,\theta_0 )$, $h \in H$ and $c > 0$ a
positive constant,
\begin{enumerate}[(1)]
\item[(1)] if $\theta_0>0$ then
%
\begin{eqnarray}
\label{Poissonupperbound} && d_{bW} \bigl(\sqrt{n}\bigl(\hat{\theta}_n(
\mathbf{X}) - \theta _0\bigr),K \bigr)\nonumber\\
&&\quad\leq\frac{2c}{\sqrt{n}} +
\frac{1}{\sqrt{n}} \biggl[2 + \frac
{(3\theta_0 + 1)^{{3}/{4}}}{\theta_0^{{3}/{4}}} \biggr]
\\
&&\qquad{}+\frac{8\theta_0}{n (\theta_0 + {c}/{n}
)^2} + \frac{\theta_0}{\sqrt{n} (\theta_0 + {c}/{n} )} + \frac
{12}{\sqrt{n}(\theta_0 + {c}/{n})}
\biggl[\frac{\theta_0}{n} + 3\theta ^2_0
\biggr]^{{1}/{2}};\nonumber
\end{eqnarray}
\item[(2)] if $\theta_0=0$ then
\[
d_{bW} \bigl(\sqrt{n}\hat{\theta}_n(\mathbf
{X}),K \bigr)= 0.
\]
\end{enumerate}
\end{Corollary}
%
\begin{rem}
(1) The upper bound expressed in \eqref{Poissonupperbound} for
the distributional distance between the actual distribution of the MLE
and the normal distribution in the case of i.i.d. random variables
following the $\operatorname{Poisson}(\theta)$ distribution, with $\theta\in
[0,\infty
)$ is of order at most $\frac{1}{\sqrt{n}}$.

(2) Since the MLE is unique and equal to $\hat{\theta
}_n(\mathbf{X}) = \bar{X}$, Lemma~\ref{Gesinetheorem} could be used
directly for~$\bar{X}$. Define $W = \sqrt{n}(\bar{X} - \theta_0) =
\frac
{1}{\sqrt{n}}\sum_{i=1}^{n} Y_i$, where\vspace*{-2pt} $Y_i = X_i - \theta_0$ are
independent, zero mean random variables. Also, ${\mathrm{E}}(W) = 0$ and
$\operatorname{Var}(W) = n\operatorname{Var}(\bar{X}) = \frac
{1}{n}\sum_{i=1}^{n}\operatorname{Var}(X_i) = \theta_0$. Therefore,
\eqref{equationLemma} for $K \sim
{\mathrm{N}}(0,\theta_0)$ and H\"{o}lder's inequality give for
$\theta_0>0$
\begin{eqnarray*}
d_{bW} \bigl(\sqrt{n}\bigl(\hat{
\theta}_n(\mathbf{X}) - \theta _0\bigr),K \bigr)
\leq\frac{1}{\sqrt{n}} \biggl(2 + \frac{1}{\theta
_0^{{3}/{2}}} \bigl[{\mathrm{E}}(Y_1)^4
\bigr]^{{3}/{4}} \biggr) = \frac
{1}{\sqrt{n}} \biggl(2 + \frac{(3\theta_0+1)^{{3}/{4}}}{\theta
_0^{{3}/{4}}}
\biggr).
\end{eqnarray*}
This bound, obtained by the direct application of Stein's method, is
smaller than the bound given in Corollary~\ref{Poissonbound}. However,
the interest in the example treated in this section, where $\Theta=
[0,\infty)$, is in adapting the approach to such cases where the MLE
could be on the boundary of the parameter space with positive
probability when it is not assumed that the MLE is a sum of random variables.
\end{rem}
%
\section{Bounds on the Mean Squared Error of the MLE}
\label{sec:Implicit}
This section focuses on the situation when an analytic form for the MLE
is not available. In the proof for the final upper bound in Theorem~\ref
{Theoremnoncan}, an explicit form of the MLE was not used. However, if
the MLE is not known, then the MSE, ${\mathrm{E}} (\hat{\theta
}_n(\mathbf{X}) - \theta_0 )^2$, appearing in the bound for~\eqref{interest} should be bounded by a quantity which is independent
of $\hat{\theta}_n(\mathbf{X})$.

Let $X_1, X_2, \ldots, X_n$ be i.i.d. random variables. Apart from the
regularity conditions, first defined in Section~\ref{sec:intro}, we
make the following further assumptions that make the steps and the
calculations easier and ensure a meaningful upper bound:
\begin{enumerate}[(Fur.1)]
\item[(Fur.1)] The support, $S$, is bounded;
\item[(Fur.2)] For $\varepsilon= \varepsilon(\theta_0)>0$ such that
$(\theta
_0-\varepsilon,\theta_0+\varepsilon)\subset\Theta$, we require
that there is a constant $C_1 = C_1(\theta_0)$ which depends on
the unknown parameter $\theta_0$ such that
 ${\sup
}_{\theta:|\theta-\theta_0|\leq\varepsilon}\llvert  l^{(3)}(\theta
; x_1)\rrvert  \leq C_1$, where $C_1 = C_1(\theta
_0)$ is a constant that depends on the unknown parameter $\theta_0$;
\item[(Fur.3)] $\exists N \in\mathbb{N}$ such that $\forall n \geq N$
we have $1 - 2\frac{\|x^2\|}{n i(\theta_0)\varepsilon^2} - \frac{\|
x\|
C_1}{\sqrt{n}[i(\theta_0)]^{{3}/{2}}} > 0$ for $\varepsilon$ as in
(Fur.2). Solving the quadratic inequality, with unknown the $\sqrt{n}$
yields that $n$, the sample size, should satisfy
\[
n \geq\frac{\|x\|^2 [C_1\varepsilon+ \sqrt
{(C_1\varepsilon)^2
+ 8[i(\theta_0)]^2} ]^2}{4[i(\theta_0)]^3\varepsilon^2}.
\]
\end{enumerate}
For ease of presentation, let $D_1 = D_1(\theta_0,x,n) = 1 - 2\frac{\|
x^2\|}{n i(\theta_0)\varepsilon^2} - \frac{\|x\|C_1}{\sqrt
{n}[i(\theta
_0)]^{{3}/{2}}}$.
%
\begin{thmm}
\label{implicitMLEbound}
Let $X_1, X_2, \ldots, X_n$ be i.i.d. random variables with density or
frequency function $f(x_i|\theta)$. Assume that the regularity
conditions \textup{(R1)--(R4)}, as well as the assumptions \textup{(Fur.1)--(Fur.3)}
are satisfied. Also assume that the MLE exists and that it is unique.
Then $A_1 = A_1(\theta_0,n)$ is an upper bound for $\sqrt{{\mathrm
{E}}(\hat
{\theta}_n(\mathbf{X}) - \theta_0)^2}$, where for $\varepsilon$
as in \textup{(Fur.2)},
%
\begin{eqnarray}
\label{A1beta} A_1 &=& [2D_1 ]^{-1} \biggl\{
\frac{2\|x\|\sqrt
{\operatorname{Var}[l''(\theta_0;X_1)]}}{n[i(\theta_0)]^{{3}/{2}}}
\nonumber
\\
& &{} + \biggl[4\frac{\|x\|^2\operatorname{
Var[l''(\theta
_0;X_1)]}}{n^2[i(\theta_0)]^3}\\
&&{}+\frac{4D_1}{n i(\theta_0)} \biggl[1+2
\frac
{\|x\|}{\sqrt{n}} \biggl(2 + \frac{{\mathrm{E}}|l'(\theta
_0;X_1)|^3}{[i(\theta
_0)]^{{3}/{2}}} \biggr) \biggr]
\biggr]^{{1}/{2}} \biggr\}.\nonumber
\end{eqnarray}
\end{thmm}
\begin{pf}
Using the notations for the remainder terms, the triangle inequality,
conditional expectations, Markov's inequality and Stein's method, the
same way as in Section~\ref{sec:generalc}, gives
\begin{eqnarray*}
&&\bigl\llvert {\mathrm{E}} \bigl[ h \bigl(\bigl(\hat{\theta
}_n(\mathbf{X}) - \theta_0\bigr)\sqrt{n i(
\theta_0)} \bigr) \bigr] - {\mathrm {E}}\bigl[h(Z)\bigr]\bigr\rrvert
\\
&&\quad\leq\frac{\|h'\|}{\sqrt{n}} \biggl(2 + \frac{{\mathrm{E}}|l'(\theta
_0;X_1)|^3}{[i(\theta_0)]^{{3}/{2}}} \biggr)
\\
&&\qquad{}+ \frac{2\|h\|{\mathrm{E}} (\hat{\theta
}_n(\mathbf{X}) -
\theta_0 )^2}{\varepsilon^2}\\
&&\qquad{}+ \frac{\|h'\|}{\sqrt{n i(\theta
_0)}}\bigl\llvert {\mathrm{E}}
\bigl[R_2(\theta_0;\mathbf{X})| \bigl|\hat {\theta
}_n(\mathbf{X}) - \theta_0\bigr|\leq\varepsilon\bigr]
\bigr\rrvert \Prob \bigl(\bigl|\hat {\theta}_n(\mathbf{X}) -
\theta_0\bigr|\leq\varepsilon \bigr)
\\
&&\qquad{} + \frac{\|h'\|}{2\sqrt{n i(\theta_0)}}{\mathrm{E}} \Bigl(\bigl(\hat {
\theta}_n(\mathbf{X}) - \theta_0\bigr)^2
\sup_{\theta:|\theta
-\theta
_0|\leq\varepsilon}\bigl\llvert l^{(3)}(\theta;\mathbf{X})
\bigr\rrvert \big| \bigl|\hat{\theta}_n(\mathbf{X}) -
\theta_0\bigr|\leq\varepsilon \Bigr).
\end{eqnarray*}
Using the definition of $R_2(\theta_0;\mathbf{x})$ and the Cauchy--Schwarz inequality yields
\begin{eqnarray*}
&&\bigl\llvert {\mathrm{E}}\bigl[R_2(\theta_0;
\mathbf{X})| \bigl|\hat {\theta}_n(\mathbf{X}) -
\theta_0\bigr|\leq\varepsilon\bigr]\bigr\rrvert \Prob \bigl(\bigl|\hat{
\theta}_n(\mathbf{X}) - \theta_0\bigr|\leq\varepsilon
\bigr)\\
&&\quad\leq {\mathrm{E}}\bigl\llvert \bigl(n i(\theta_0) +
l''(\theta_0;\mathbf {X})\bigr) \bigl(
\hat{\theta }_n(\mathbf{X}) - \theta_0\bigr)\bigr
\rrvert
\\
&&\quad \leq\sqrt{{\mathrm{E}}\bigl[n i(\theta_0) +
l''(\theta _0;\mathbf {X})
\bigr]^2{\mathrm{E}}\bigl[\hat{\theta}_n(\mathbf{X})
- \theta_0\bigr]^2}\\
&&\quad = \sqrt {n\operatorname{Var}
\bigl(l''(\theta_0;X_1)\bigr)}
\sqrt{{\mathrm{E}}\bigl[\hat {\theta}_n(\mathbf {X})-
\theta_0\bigr]^2},
\end{eqnarray*}
which leads to
%
\begin{eqnarray}
\label{boundbetaim}
& &\bigl\llvert {\mathrm{E}} \bigl[ h \bigl(\bigl(\hat{
\theta }_n(\mathbf{X}) - \theta_0\bigr)\sqrt{n i(
\theta_0)} \bigr) \bigr] - {\mathrm {E}}\bigl[h(Z)\bigr]\bigr\rrvert
\nonumber\\
&&\quad\leq\frac{\|h'\|}{\sqrt{n}} \biggl(2 + \frac{{\mathrm{E}}|l'(\theta
_0;X_1)|^3}{[i(\theta_0)]^{{3}/{2}}} \biggr)
\nonumber
\\[-8pt]
\\[-8pt]
\nonumber
&&\qquad{}+ \frac{2\|h\|{\mathrm{E}} (\hat{\theta}_n(\mathbf{X}) -
\theta
_0 )^2}{\varepsilon^2} + \frac{\|h'\|nC_1}{2\sqrt{n i(\theta
_0)}}{\mathrm{E}} \bigl(\hat{
\theta}_n(\mathbf{X}) - \theta _0
\bigr)^2
\\
&&\qquad{}+ \frac{\|h'\|\sqrt{\operatorname{Var}(l''(\theta
_0;X_1))}\sqrt{{\mathrm{E}} (\hat{\theta}_n(\mathbf{X}) -
\theta_0 )^2}}{\sqrt
{i(\theta_0)}}.\nonumber
\end{eqnarray}
Straightforward calculations and denoting with $B_{x^2}$ the upper
bound for \eqref{interest} when $h(x) = x^2$, lead to
%
\begin{eqnarray}
\label{quadineq}
{\mathrm{E}} \bigl(\hat{\theta}_n(
\mathbf{X}) - \theta_0 \bigr)^2 &=&
\frac{1}{n i(\theta_0)}\bigl\llvert {\mathrm{E}} \bigl[\sqrt{n i(\theta
_0)}\bigl(\hat{\theta}_n(\mathbf{X}) -
\theta_0\bigr) \bigr]^2 - {\mathrm{E}}
\bigl(Z^2 \bigr) + {\mathrm{E}} \bigl(Z^2 \bigr)\bigr
\rrvert
\nonumber
\\[-8pt]
\\[-8pt]
\nonumber
& \leq&\frac{1}{n i(\theta_0)}(B_{x^2} + 1),
\end{eqnarray}
where
\begin{eqnarray*}
B_{x^2} &\leq& 2\frac{\|x\|}{\sqrt{n}} \biggl(2 +
\frac
{{\mathrm{E}}|l'(\theta_0;X_1)|^3}{[i(\theta_0)]^{{3}/{2}}} \biggr) + \frac{2\|
x^2\|{\mathrm{E}} (\hat{\theta}_n(\mathbf{X}) - \theta
_0
)^2}{\varepsilon^2} + \frac{\|x\|\sqrt{n}C_1}{\sqrt{i(\theta
_0)}}{
\mathrm{E}} \bigl(\hat{\theta}_n(\mathbf{X}) - \theta
_0 \bigr)^2
\\
&&{} + 2\frac{\|x\|\sqrt{\operatorname{Var}(l''(\theta
_0;X_1))}\sqrt{{\mathrm{E}} (\hat{\theta}_n(\mathbf{X}) -
\theta_0 )^2}}{\sqrt
{i(\theta_0)}}.
\end{eqnarray*}
Now $B_{x^2}$ also includes ${\mathrm{E}} (\hat{\theta
}_n(\mathbf
{X}) - \theta_0 )^2$ and its positive root. Therefore, the next
step is to solve the simple quadratic inequality \eqref{quadineq}, with
unknown $\sqrt{{\mathrm{E}} (\hat{\theta}_n(\mathbf{X}) -
\theta
_0 )^2}$. Using (Fur.3), after basic calculations we obtain that
$0 < \sqrt{{\mathrm{E}} (\hat{\theta}_n(\mathbf{X}) -
\theta_0
)^2} \leq A_1$.
\end{pf}
%
\begin{rem} (1) Using this result, the final upper bound for
\eqref{interest} which is useful when no analytic expression of the MLE
is available, becomes\vspace*{-2pt}
\begin{eqnarray}
\label{implicit}
\nonumber
d_{bW} \bigl(\sqrt{n i(\theta_0)}
\bigl(\hat{\theta }_n(\mathbf {X}) - \theta_0
\bigr),Z \bigr) &\leq&\frac{1}{\sqrt{n}} \biggl(2 + \frac
{{\mathrm{E}}|l'(\theta_0;X_1)|^3}{[i(\theta_0)]^{{3}/{2}}} \biggr) +
\frac
{2(A_1)^2}{\varepsilon^2}
\nonumber
\\[-8pt]
\\[-8pt]
\nonumber
&&{} + \frac{\sqrt{n}C_1(A_1)^2}{2\sqrt{i(\theta_0)}} + \frac{\sqrt
{\operatorname{Var}[l''(\theta_0;X_1)]}A_1}{\sqrt{i(\theta_0)}}.
\end{eqnarray}

(2) The\vspace*{-2pt} order of $A_1$ in terms of the sample size is $\frac
{1}{\sqrt{n}}$ and hence the order of the final upper bound in \eqref
{implicit} is also $\frac{1}{\sqrt{n}}$.
\end{rem}
\begin{example*}[(The Beta distribution)]
Consider the example of i.i.d random variables from the Beta
distribution with one of the two shape parameters being unknown. In
this case, the MLE can only be expressed in terms of the inverse of the
digamma function, $\Psi(\theta) = \frac{\mathrm{d}}{\mathrm
{d}\theta
}\log\Gamma(\theta)$. We use the general result in Theorem~\ref
{implicitMLEbound}, in order to obtain an upper bound for the MSE and
use it to get an upper bound for \eqref{interest}. The following
corollary gives the result.
\end{example*}
%
\begin{Corollary}
\label{Corollary_BETA}
Let $X_1, X_2, \ldots, X_n$ be i.i.d. random variables from
the $\operatorname{Beta}(\theta_0, \beta)$ distribution, where $\beta$ is
known and $\theta_0$ is unknown. Let $B_1 = B_1(\theta_0) =
8 (\Psi_3(\theta_0) + \Psi_3(\theta_0+\beta) + 3[\Psi
_1(\theta
_0)]^2 + 3[\Psi_1(\theta_0+\beta)]^2 )$, where $\Psi
_j(\theta), j\in\mathbb{N}$ is the $j$th derivative of the digamma
function, $\Psi(\theta)$. Also, let $B_2=B_2(\theta_0) = \frac
{96\beta+ 6.6\beta\theta_0^4}{\theta_0^4}$, $D_{\Psi1}=D_{\Psi
1}(\theta_0,\beta) = \Psi_1(\theta_0) - \Psi_1(\theta_0 + \beta
)$ and\vspace*{-2pt}
%
\begin{eqnarray}
\label{B3beta} B_3 = B_3(\theta_0,n)& =&
\biggl[ \biggl(4+\frac{8}{\sqrt
{n}} \biggl(2 +
\frac{(B_1)^{{3}/{4}}}{D_{\Psi1}^{{3}/{2}}} \biggr)
\biggr)
\biggl(1-\frac{8}{n\theta_0^2D_{\Psi1}}-\frac{B_2}{\sqrt{n}D_{\Psi
1}^{{3}/{2}}}  \biggr)  \biggr]^{{1}/{2}}
\nonumber
\\[-9pt]
\\[-9pt]
\nonumber
&&{}\times\biggl(2 \biggl(\sqrt{D_{\Psi1}} -
\frac
{8}{n\theta_0^2\sqrt{D_{\Psi1}}} - \frac{B_2}{\sqrt{n}D_{\Psi
1}} \biggr)\biggr)^{-1}.\vspace*{-1pt}
\end{eqnarray}
Let\vspace*{-2pt}
\begin{eqnarray*}
n \geq\biggl[B_2\frac{\theta_0}{2} + \sqrt{\frac{(B_2\theta
_0)^2}{4} + 8\bigl[\Psi_1(\theta_0) - \Psi_1(\theta_0 + \beta
)\bigr]^2}
\biggr]^2\bigl(\bigl[\Psi_1(\theta_0) - \Psi_1(\theta_0 +
 \beta)\bigr]^3\theta_0^2\bigr)^{-1}.
\end{eqnarray*}
Then for $Z \sim \mathrm{N}(0,1)$\vspace*{-1pt}
%
\begin{eqnarray}
\label{boundbetaim2} && d_{bW} \bigl(\sqrt{n i(\theta_0)}\bigl(
\hat{\theta}_n(\mathbf{X}) - \theta_0\bigr),Z \bigr)
\leq\frac{1}{\sqrt{n}} \biggl(2 + \frac
{(B_1)^{{3}/{4}}}{[\Psi_1(\theta_0) - \Psi_1(\theta_0+\beta)]^{{3}/{2}}} \biggr)\nonumber
\\[-9pt]
\\[-9pt]
&&\qquad{} + \frac{8}{n\theta_0^2}(B_3)^2 +
\frac
{B_2(B_3)^2}{2\sqrt
{n}[\Psi_1(\theta_0) - \Psi_1(\theta_0+\beta)]^{{1}/{2}}}.\nonumber
\end{eqnarray}
\end{Corollary}
\begin{pf}
See the \hyperref[app]{Appendix}.
\end{pf}

Now, we study the accuracy of our bound for the MSE of the MLE by
simulations. For the simulations, $\theta_0=1.5$, $\beta=1$ and in this
case of $\beta$ being equal to 1, the MLE is $\hat{\theta
}_n(\mathbf
{X}) = - \frac{n}{\sum_{i=1}^{n}\log X_i}$. We find that $n \geq7460$,
in order for (Fur.3) to be satisfied. The process to simulate is quite
simple. Let $n \in \{7460, 7461, \ldots, 8459  \}$
and for each $n$, start by generating 10\,000 trials of $n$ random
independent observations, $x$, from the Beta distribution with
parameter values as above. We evaluate the MLE, $\hat{\theta
}_n(\mathbf{X})$, of the parameter in each trial, which in turn
gives a vector of 10\,000 values. Thus, for each $n$ from 7460 to 8459,
we evaluate the sample MSE, $\hat\mathrm{{E}} (\hat{\theta
}_n(\mathbf{X}) - \theta_0 )^2 = \frac{1}{10\,000}\sum_{i=1}^{10\,000} [\hat{\theta}_n(\mathbf{x})[i] - \theta
_0
]^2$ and compare it with its upper bound, $ (\frac{B_3}{\sqrt
{n}} )^2$, where $B_3$ is given in \eqref{B3beta}. The difference
between their values measures the error of our bound on the MSE. Part
of the results from the simulations is shown in Table~\ref{Tablebeta}.
The table indicates that the bound and the error decrease as the sample
size increases, as expected, since the order of the upper bound for the
MSE is $\frac{1}{n}$. In addition, it is reasonable that the smaller
the sample size is, the larger the bound is. The bounds are
considerably larger than the estimated MSE and they are not numerically
sharp. In addition, because of the relatively strong requirement that
$n\geq7460$, these bounds on the MSE are more of theoretical
interest.

%
\begin{table}
\caption{Part of the results taken by simulations from the
$\operatorname{Beta}(1.5,1)$ distribution}\label{Tablebeta}
\begin{tabular*}{\textwidth}{@{\extracolsep{\fill}}llll@{}}
\hline
\multicolumn{1}{@{}l}{$n$} & $\hat\mathrm{{E}} (\hat{\theta}_n(\mathbf{X}) -
\theta
_0 )^2$ & Upper bound & Error\\
\hline
7500 & 0.0002 & 0.2517 & 0.2515\\
7700 & 0.0002 & 0.0416 & 0.0414\\
7900 & 0.0002 & 0.0223 & 0.0221\\
8100 & 0.0002 & 0.0151 & 0.0149\\
8300 & 0.0002 & 0.0112 & 0.00110\\
\hline
\end{tabular*}
\end{table}

\begin{remarks*}
Several interesting
paths lead from the work explained in this paper. When the dimension of
the parameter is $d>1$, Stein bounds are available in Chen
\textit{et~al.} \cite{Chen_book},
which can be employed to get upper bounds related to the distribution
of the MLE in a multi-parameter setting (work in progress). In
addition, one of the main advantages of Stein's method is that it can
be used in situations where dependence comes into play. Upper bounds on
the distributional distance between the distribution of the MLE and the
normal distribution in the case of dependent random variables are also
work in progress.
\end{remarks*}

\begin{appendix}\label{app}
\section*{Appendix: Some proofs}
\begin{pf*}{Proof of Lemma~\ref{Lemmaincreasing}}
Let $\varepsilon> 0$ and $f$ a continuous
increasing function with $f(m) \geq0$ for $m > 0$. Then,
\begin{eqnarray*}\label{conditional}
{\mathrm{E}} \bigl[f(M) \bigr] &= &{
\mathrm{E}} \bigl[f(M)| M \leq \varepsilon \bigr]\Prob (M \leq\varepsilon )
+ {\mathrm {E}} \bigl[f(M)|M>\varepsilon \bigr]\Prob (M > \varepsilon )
\\
&=& {\mathrm{E}} \bigl[f(M)| M \leq\varepsilon \bigr] \bigl(1 -
\Prob (M > \varepsilon ) \bigr) + {\mathrm{E}} \bigl[f(M) |M>\varepsilon
\bigr]\Prob (M > \varepsilon )
\\
& =& {\mathrm{E}} \bigl[f(M)| M \leq\varepsilon \bigr] + \Prob (M
> \varepsilon ) \bigl({\mathrm{E}} \bigl[f(M) |M>\varepsilon \bigr] - {
\mathrm{E}} \bigl[f(M)| M \leq\varepsilon \bigr] \bigr)
\\
& \geq&{\mathrm{E}} \bigl[f(M)| M \leq\varepsilon \bigr]\qquad \mbox{as }f(m) \mbox{ is increasing.}
\end{eqnarray*}
\upqed\end{pf*}

\begin{pf*}{Proof of Corollary~\ref{Corollary_BETA}} The probability density function is
%
\setcounter{equation}{0}
\begin{equation}
\label{densitybeta} f(x|\theta) = \frac{\Gamma(\theta+\beta)}{\Gamma(\theta)\Gamma
(\beta
)}x^{\theta-1}(1-x)^{\beta-1},
\end{equation}
with $\theta> 0$ and $x \in[0,1]$. Hence
%
\begin{eqnarray}
\label{betalikelihood} l(\theta;\mathbf{x}) &= &n\bigl[\log\bigl(\Gamma(\theta+
\beta)\bigr) - \log \bigl(\Gamma (\theta)\bigr) - \log\bigl(\Gamma(\beta)\bigr)
\bigr]
\nonumber
\\[-8pt]
\\[-8pt]
\nonumber
&&{}+ (\theta- 1)\sum_{i=1}^{n}\log
x_i +(\beta-1)\sum_{i=1}^{n}
\log(1-x_i)
\end{eqnarray}
and
\begin{eqnarray*}
l'(\theta;\mathbf{x}) &=& n\bigl[\Psi(\theta+\beta) - \Psi(
\theta)\bigr] + \sum_{i=1}^{n} \log
x_i
\\
l^{(j)}(\theta;\mathbf{x})& =& n\bigl(\Psi_{j-1}(\theta+
\beta) - \Psi _{j-1}(\theta)\bigr),\qquad  j\in\mathbb{N}\setminus\{1\}.
\end{eqnarray*}
Now we show that the conditions (R1)--(R4) and the assumptions
(Fur.1)--(Fur.3) are satisfied. For (R1) it is obvious. As for (R2), the
three times differentiability of the density function can be verified
from \eqref{betalikelihood}. In addition, using \eqref{densitybeta} and
the expressions for the logarithmic expectations of a Beta distributed
random variable, it is straightforward\vspace*{-1pt} to verify $\bigintsss
_{0}^{1}\frac{\mathrm{d}^j}{\mathrm{d}\theta^j} f(x|\theta)
\,\mathrm
{d}x = \frac{\mathrm{d}^j}{\mathrm{d}\theta^j}\bigintsss_{0}^{1}
f(x|\theta) \,\mathrm{d}x = 0, j \in \{1, 2, 3 \}$
for (R2). Let $\varepsilon= \varepsilon(\theta_0) > 0$ such that
$\theta\in
(\theta_0 - \varepsilon, \theta_0 + \varepsilon) \subset\Theta$.
Since in
this case $\Theta= (0,\infty)$, indeed $0 < \varepsilon< \theta_0$.
Using a first order Taylor expansion and the fact that
%
\begin{equation}
\label{psi_m} \Psi_m(z) = (-1)^{m+1}m!\sum
_{k=0}^{\infty}\frac{1}{(z+k)^{m+1}}
\end{equation}
gives
\[
\Psi_3(z) = 6\sum_{k=0}^{\infty}
\frac{1}{(z+k)^4}\qquad \mbox{for } z \in \mathbb{C}\setminus\bigl\{\mathbb{Z^-}\bigr
\}\mbox{ and } m>0,
\]
with $\Psi_3(z)$ being a decreasing function of $z$. For $\theta\in
(\theta_0 - \varepsilon, \theta_0 + \varepsilon)$,
%
\begin{eqnarray}
\label{R4beta} \biggl\llvert \frac{\mathrm{d}^3}{\mathrm{d}\theta^3}\log f(x|\theta )\biggr\rrvert &=&
\bigl\llvert \Psi_2(\theta+ \beta) - \Psi_2(\theta)\bigr
\rrvert
\nonumber
\\[-8pt]
\\[-8pt]
\nonumber
& \leq&\beta\bigl\llvert \Psi_3\bigl(\theta^*\bigr)\bigr
\rrvert \leq\beta \bigl\llvert \Psi_3(\theta_0-
\varepsilon)\bigr\rrvert = M(x),
\end{eqnarray}
with ${\mathrm{E}}[M(X)] < \infty$. Hence, (R3) holds as well. Also,
$i(\theta_0) = \Psi_1(\theta_0) - \Psi_1(\theta_0 + \beta)$ which is
positive since it is obvious from \eqref{psi_m} that $\Psi_1(z)$ is a
decreasing function. The assumption (Fur.1) obviously holds with $\|x\|
\leq1$. Using \eqref{R4beta} and the fact that $\sum_{i=1}^{\infty
}\frac{1}{i^4} = \frac{\pi^4}{90} < 1.1$ gives
%
\begin{eqnarray}
\label{C1beta}
\nonumber
\sup_{\theta:|\theta-\theta_0|\leq\varepsilon}\bigl\llvert l^{(3)}(
\theta;X_1)\bigr\rrvert &\leq&\beta\bigl\llvert \Psi_3(
\theta _0-\varepsilon )\bigr\rrvert
= 6\beta\sum_{k=0}^{\infty}
\frac{1}{(\theta_0 -
\varepsilon+
k)^4}
\nonumber
\\[-8pt]
\\[-8pt]
\nonumber
&\leq&6\beta \Biggl[\frac{1}{(\theta_0-\varepsilon)^4} + \sum
_{k=1}^{\infty}\frac{1}{k^4} \Biggr]
<\frac{6\beta}{(\theta_0 - \varepsilon)^4} + 6.6\beta= C_1.\nonumber
\end{eqnarray}
Thus, (Fur.2) is also satisfied. Now, since $i(\theta_0) = \Psi
_1(\theta
_0) - \Psi_1(\theta_0 + \beta)$ take
\[
n \geq\frac{ [C_1\varepsilon+ \sqrt{(C_1\varepsilon
)^2 +
8[\Psi_1(\theta_0) - \Psi_1(\theta_0 + \beta)]^2}
]^2}{4\varepsilon
^2[\Psi_1(\theta_0) - \Psi_1(\theta_0 + \beta)]^3}
\]
in order for (Fur.3) to be satisfied. To find $B_3$, firstly, as
${\mathrm{E}}|l'(\theta_0;X_1)|^3$ is not straightforward to evaluate
due to the
absolute value in the expectation, it is easily seen that using H\"
{o}lder's inequality ${\mathrm{E}}|l'(\theta_0;X_1)|^3 \leq
[{\mathrm{E}}(l'(\theta_0;X_1))^4 ]^{{3}/{4}}$ we find an
upper bound for
\begin{eqnarray*}
{\mathrm{E}}\bigl[l'(\theta_0;X_1)
\bigr]^4 &=& {\mathrm{E}}\bigl[\log X_1 + \Psi (
\theta_0 + \beta) - \Psi(\theta_0)\bigr]^4 \\
&=&
{\mathrm{E}}\bigl[\log X_1 - {\mathrm{E}}(\log X_1)
\bigr]^4.
\end{eqnarray*}
If $G_1 \sim\Gamma(\theta_0,\lambda)$ and $G_2 \sim\Gamma(\beta
,\lambda)$ independent, then $\frac{G_1}{G_1 + G_2} \sim
\operatorname{
Beta}(\theta_0,\beta)$. Thus, with $X_1 = \frac{G_1}{G_1 + G_2}$
%
\begin{eqnarray}
\label{4momentbeta}
\hspace*{-25pt}\nonumber {\mathrm{E}} \bigl[l'(
\theta_0;X_1) \bigr]^4
&=& {\mathrm {E}}
\bigl[\bigl(\log G_1 - {\mathrm{E}}[\log G_1]\bigr) +
\bigl({\mathrm{E}}\bigl[\log(G_1 + G_2)\bigr] -
\log(G_1 + G_2)\bigr) \bigr]^4
\\[-8pt]
\\[-8pt]
\hspace*{-25pt}&\leq& 8 \bigl[{\mathrm{E}} \bigl(\log G_1 - {\mathrm{E}}(\log
G_1) \bigr)^4 + {\mathrm{E}} \bigl(\log(G_1
+ G_2) - {\mathrm{E}}\bigl(\log(G_1 + G_2)
\bigr) \bigr)^4 \bigr].\nonumber
\end{eqnarray}
Now we calculate the fourth central moment of the logarithm of a Gamma
distributed random variable. Using that $\bigintsss_{0}^{\infty}
\frac
{z^{\alpha-1}\mathrm{ e}^{-z}(\log z)^k}{\Gamma(\alpha)} \,\mathrm
{d}z =
\frac{\Gamma^{(k)}(\alpha)}{\Gamma(\alpha)}$, for any $\alpha> 0$ and
$k \in\mathbb{N}$ gives that for \mbox{$Y \sim\Gamma(\alpha, \lambda)$}
\[
{\mathrm{E}}(\log Y) = \Psi(\alpha) - \log\lambda.
\]
Using again $z=\lambda y$,
%
\begin{eqnarray*}
\nonumber
{\mathrm{E}} \bigl[\log Y - {\mathrm{E}}(\log Y) \bigr]^4
&=& \bigintsss _{0}^{\infty} \frac{z^{\alpha-1}\mathrm{ e}^{-z}}{\Gamma(\alpha
)} \biggl(\log
\biggl(\frac{z}{\lambda} \biggr) - {\mathrm{E}} \biggl(\log \biggl(
\frac
{Z}{\lambda} \biggr) \biggr) \biggr)^4 \,\mathrm{d}z
\\
& =& \bigintsss_{0}^{\infty} \frac{z^{\alpha-1}\mathrm{
e}^{-z}}{\Gamma(\alpha)} \bigl(
\log z - {\mathrm{E}} (\log Z ) \bigr)^4 \,\mathrm{d}z
\nonumber
\\[-8pt]
\\[-8pt]
\nonumber
&= &\frac{1}{\Gamma(\alpha)}\sum_{k=0}^{4}
\pmatrix{4
\cr
k}(-1)^k\bigl[\Psi(\alpha)\bigr]^{4-k}
\bigintsss_{0}^{\infty} z^{\alpha
-1}\mathrm{
e}^{-z}(\log z)^k \,\mathrm{d}z
\\
\nonumber
& =& -3\bigl[\Psi(\alpha)\bigr]^4 + 6\bigl[\Psi(\alpha)
\bigr]^2\frac{\Gamma
''(\alpha)}{\Gamma(\alpha)} - 4\Psi(\alpha)\frac{\Gamma
^{(3)}(\alpha
)}{\Gamma(\alpha)} +
\frac{\Gamma^{(4)}(\alpha)}{\Gamma(\alpha)}.
\end{eqnarray*}
At this point, the digamma function can be used in order to simplify
the expression above. Following simple steps it can be easily verified that
\begin{eqnarray*}
\frac{\Gamma''(\alpha)}{\Gamma(\alpha)} &=& \Psi _1(\alpha) + \bigl[\Psi(
\alpha)\bigr]^2,\qquad \frac{\Gamma^{(3)}(\alpha)}{\Gamma(\alpha)} = \Psi_2(\alpha) + 3
\Psi(\alpha)\Psi_1(\alpha) + \bigl[\Psi(\alpha )\bigr]^3,
\\
\frac{\Gamma^{(4)}(\alpha)}{\Gamma(\alpha)} &=& \Psi _3(\alpha ) + 4
\Psi_2(\alpha)\Psi(\alpha) + 6\Psi_1(\alpha)\bigl[\Psi(
\alpha )\bigr]^2 + 3\bigl[\Psi_1(\alpha)
\bigr]^2 + \bigl[\Psi(\alpha)\bigr]^4.
\end{eqnarray*}
Hence for $Y \sim\Gamma(\alpha,\lambda)$
\[
{\mathrm{E}} \bigl[\log Y - {\mathrm{E}}(\log Y) \bigr]^4 =
\Psi_3(\alpha ) + 3\bigl[\Psi_1(\alpha)
\bigr]^2
\]
and therefore, from \eqref{4momentbeta},
\[
{\mathrm{E}} \bigl[l'(\theta_0;X_1)
\bigr]^4 \leq8 \bigl(\Psi _3(\theta_0) +
\Psi_3(\theta_0 + \beta) + 3 \bigl[\Psi_1(
\theta _0) \bigr]^2 + 3 \bigl[\Psi_1(
\theta_0 + \beta) \bigr]^2 \bigr) = B_1.
\]
With $C_1$ as in \eqref{C1beta}, taking $\varepsilon= \frac{\theta
_0}{2}$, we conclude that
\[
\sup_{\theta:|\theta-\theta_0|\leq\varepsilon}\bigl\llvert l^{(3)}(
\theta;X_1)\bigr\rrvert \leq\frac{96\beta}{\theta_0^4} + 6.6\beta=
B_2.
\]
Using \eqref{betalikelihood}, gives
\[
\operatorname{Var}\bigl(l''(
\theta_0;X_1)\bigr) = \operatorname {Var}\bigl(
\Psi_1(\theta_0 + \beta) - \Psi_1(
\theta_0)\bigr)=0.
\]
Having found all the necessary quantities, we calculate the upper bound
in \eqref{A1beta} and multiply it by $\sqrt{n}$. This is equal to $B_3$
shown in \eqref{B3beta}, which is an upper bound for $\sqrt{n{\mathrm
{E}}(\hat{\theta}_n(\mathbf{X}) - \theta_0)^2}$ in the specific case
of i.i.d. random variables from the Beta distribution. Using this bound
in \eqref{boundbetaim} gives the result in \eqref{boundbetaim2}.
\end{pf*}
\end{appendix}

\section*{Acknowledgements}
The authors would like to thank Robert E. Gaunt for various insightful
comments and for his idea regarding a perturbation to treat the Poisson
example. This was the motivation behind the necessity to generalize the
perturbation idea in Section~\ref{sec:perturbationc} for any discrete
distribution. We would also like to thank Larry Goldstein for his
comments and Adrian R\"{o}llin for suggesting the Poisson example. In
addition, the authors sincerely thank two anonymous reviewers for
suggestions that lead to an improvement of the paper. Andreas
Anastasiou was supported by a Teaching Assistantship Bursary from the
Department of Statistics, University of Oxford, and an Engineering and
Physical Sciences Research Council (EPSRC) Scholarship. Gesine Reinert
was supported in part by EPSRC grant EP/K032402/1.



%
%


\printhistory
\end{document}